\documentclass[11pt]{amsart}
\usepackage{amscd}
\usepackage[arrow,matrix]{xy}
\usepackage{graphicx}
\usepackage{amsmath}
\usepackage{amsmath, latexsym, amssymb}
\input xypic
\numberwithin{equation}{section}
\theoremstyle{plain}
\newtheorem{lemma}{Lemma}[section]
\newtheorem{proposition}[lemma]{Proposition}
\newtheorem{theorem}[lemma]{Theorem}
\newtheorem{corollary}[lemma]{Corollary}

\theoremstyle{definition}
\newtheorem{definition}[lemma]{Definition}
\newtheorem{remark}[lemma]{Remark}

\begin{document}
\newcommand{\R}{{\mathbb R}}
\newcommand{\C}{{\mathbb C}}
\newcommand{\F}{{\mathbb F}}
\renewcommand{\O}{{\mathbb O}}
\newcommand{\Z}{{\mathbb Z}} 
\newcommand{\N}{{\mathbb N}}
\newcommand{\Q}{{\mathbb Q}}
\renewcommand{\H}{{\mathbb H}}

\newcommand{\Aa}{{\mathcal A}}
\newcommand{\Bb}{{\mathcal B}}
\newcommand{\Cc}{{\mathcal C}}    
\newcommand{\Dd}{{\mathcal D}}
\newcommand{\Ee}{{\mathcal E}}
\newcommand{\Ff}{{\mathcal F}}
\newcommand{\Gg}{{\mathcal G}}    
\newcommand{\Hh}{{\mathcal H}}
\newcommand{\Kk}{{\mathcal K}}
\newcommand{\Ii}{{\mathcal I}}
\newcommand{\Jj}{{\mathcal J}}
\newcommand{\Ll}{{\mathcal L}}    
\newcommand{\Mm}{{\mathcal M}}    
\newcommand{\Nn}{{\mathcal N}}
\newcommand{\Oo}{{\mathcal O}}
\newcommand{\Pp}{{\mathcal P}}
\newcommand{\Qq}{{\mathcal Q}}
\newcommand{\Rr}{{\mathcal R}}
\newcommand{\Ss}{{\mathcal S}}
\newcommand{\Tt}{{\mathcal T}}
\newcommand{\Uu}{{\mathcal U}}
\newcommand{\Vv}{{\mathcal V}}
\newcommand{\Ww}{{\mathcal W}}
\newcommand{\Xx}{{\mathcal X}}
\newcommand{\Yy}{{\mathcal Y}}
\newcommand{\Zz}{{\mathcal Z}}

\newcommand{\zt}{{\tilde z}}
\newcommand{\xt}{{\tilde x}}
\newcommand{\Ht}{\widetilde{H}}
\newcommand{\ut}{{\tilde u}}
\newcommand{\Mt}{{\widetilde M}}
\newcommand{\Llt}{{\widetilde{\mathcal L}}}
\newcommand{\yt}{{\tilde y}}
\newcommand{\vt}{{\tilde v}}
\newcommand{\Ppt}{{\widetilde{\mathcal P}}}
\newcommand{\bp }{{\bar \partial}} 

\newcommand{\Remark}{{\it Remark}}
\newcommand{\Proof}{{\it Proof}}
\newcommand{\ad}{{\rm ad}}
\newcommand{\Om}{{\Omega}}
\newcommand{\om}{{\omega}}
\newcommand{\eps}{{\varepsilon}}
\newcommand{\Di}{{\rm Diff}}
\newcommand{\im}{{\rm Im}}
\newcommand{\Pro}[1]{\noindent {\bf Proposition #1}}
\newcommand{\Thm}[1]{\noindent {\bf Theorem #1}}
\newcommand{\Lem}[1]{\noindent {\bf Lemma #1 }}
\newcommand{\An}[1]{\noindent {\bf Anmerkung #1}}
\newcommand{\Kor}[1]{\noindent {\bf Korollar #1}}
\newcommand{\Satz}[1]{\noindent {\bf Satz #1}}

\renewcommand{\a}{{\mathfrak a}}
\renewcommand{\b}{{\mathfrak b}}
\newcommand{\e}{{\mathfrak e}}
\renewcommand{\k}{{\mathfrak k}}
\newcommand{\pg}{{\mathfrak p}}
\newcommand{\g}{{\mathfrak g}}
\newcommand{\gl}{{\mathfrak gl}}
\newcommand{\h}{{\mathfrak h}}
\renewcommand{\l}{{\mathfrak l}}
\newcommand{\sm}{{\mathfrak m}}
\newcommand{\n}{{\mathfrak n}}
\newcommand{\s}{{\mathfrak s}}
\renewcommand{\o}{{\mathfrak o}}
\newcommand{\so}{{\mathfrak so}}
\renewcommand{\u}{{\mathfrak u}}
\newcommand{\su}{{\mathfrak su}}
\newcommand{\ssl}{{\mathfrak sl}}
\newcommand{\ssp}{{\mathfrak sp}}
\renewcommand{\t}{{\mathfrak t }}
\newcommand{\Cinf}{C^{\infty}}
\newcommand{\la}{\langle}
\newcommand{\ra}{\rangle}
\newcommand{\half}{\scriptstyle\frac{1}{2}}
\newcommand{\p}{{\partial}}
\newcommand{\notsub}{\not\subset}
\newcommand{\iI}{{I}}               
\newcommand{\bI}{{\partial I}}      
\newcommand{\LRA}{\Longrightarrow}
\newcommand{\LLA}{\Longleftarrow}
\newcommand{\lra}{\longrightarrow}
\newcommand{\LLR}{\Longleftrightarrow}
\newcommand{\lla}{\longleftarrow}
\newcommand{\INTO}{\hookrightarrow}

\newcommand{\QED}{\hfill$\Box$\medskip}
\newcommand{\UuU}{\Upsilon _{\delta}(H_0) \times \Uu _{\delta} (J_0)}
\newcommand{\bm}{\boldmath}

\title[ Cohomology theories  on l.c.s.  manifolds]{\large   Cohomology theories on locally conformal  symplectic manifolds}
\author{H\^ong V\^an L\^e  and Ji\v ri Van\v zura }
\thanks{H.V.L. and J.V. are partially supported by RVO: 67985840} \maketitle

\medskip

\abstract In this note  we     introduce  primitive cohomology  groups of  locally conformal symplectic manifolds
$(M^{2n}, \om, \theta)$. We study the relation between the  primitive cohomology groups and the   Lichnerowicz-Novikov  cohomology  groups  of $(M^{2n}, \om, \theta)$, using and  extending the technique of spectral sequences developed by Di Pietro and Vinogradov for symplectic manifolds.    We  discuss  related results  by many peoples, e.g. Bouche, Lychagin, Rumin, Tseng-Yau, in light of our spectral sequences. We calculate the primitive cohomology groups  of a $(2n+2)$-dimensional locally conformal symplectic nilmanifold  as well as those of a l.c.s.  solvmanifold. We show that the  l.c.s. solvmanifold  is a mapping torus  of a contactomorphism, which is not isotopic to the identity.
\endabstract

\tableofcontents

\section{Introduction}

A differentiable manifold $(M^{2n}, \om, \theta)$ provided with a non-degenerate  2-form $\om$ and a closed 1-form $\theta$  is called a locally conformal symplectic  (l.c.s.) manifold, if   $d\om = -\om \wedge \theta$, $d\theta = 0$. The 1-form $\theta$ is called  the Lee form of $\om$.    

The class of l.c.s. manifolds has  attracted 
strong interests among geometers in recent years. For instance, Vaisman  showed that l.c.s. manifolds
may be viewed as phase spaces for a natural generalization of Hamiltonian dynamics \cite{Vaisman1985}. Bande and Kotschick showed that  a pair composed  of a contact manifold and a contactomorphism  is naturally  associated with a l.c.s. manifold \cite{BK2010} (see also Proposition \ref{mtc} and Proposition \ref{inv2} below). Furthermore, l.c.s. manifolds together with contact manifolds are the only transitive Jacobi manifolds \cite[Remark 2.10]{Marle1989}. It is also worth mentioning that locally conformal K\"ahler manifolds, a natural subclass of l.c.s. manifolds,  are actively studied in complex geometry,  e.g. see \cite{HK2011}, \cite{OV2010}.

Note that a l.c.s. manifold is locally conformal equivalent to a symplectic manifold, i.e.  locally $\theta = df$ and $ \om = e^{-f} \om_0$, $d\om_0 = 0$. By the Darboux theorem all symplectic manifolds of the same dimension are locally  equivalent. Hence  symplectic manifolds have only global invariants,  and cohomological invariants are most natural among them. First  (co)homological symplectic invariants  were proposed  in  works by Gromov and Floer then followed  by works by McDuff,  Hofer and Salamon, Fukaya and Ono,  Ruan, Tian, Witten    and many others, including the first author of this note.  This approach  was based on the use of the  theory of  elliptic  differential operators with purpose  to make regular  certain  Morse (co)homology theory or the intersection theory on the infinite dimensional loop space  on  a symplectic  manifold $M^{2n}$, or on the space of holomorphic curves on  $M^{2n}$.  This  elliptic (co)homology  theory has  huge success, but the  computational part of the theory is quite complicated.  Almost at the same time, a ``linear"  symplectic cohomology theory  has been  developed, beginning with the paper by Brylinski \cite{Brylinski1988}, followed by Bouche  \cite{Bouche1990}, and  then by  other peoples (see \cite{FIL1998}, \cite{PV2006}, \cite{TY2009}). This theory is mostly motivated by
analogues in K\"ahler geometry, the Dolbeault theory, and  the cohomology theory for differential equations developed by Vinogradov and his school.

This linear symplectic cohomology theory has not yet drawn as much attention as it, to our opinion, should have. This is, probably, due to the fact that its potentially important applications are still in a phase of elaboration. The computational part of the linear theory seems to be not so complicated as in the elliptic theory, and this is an advantage of its.

In our paper we  further develop the linear symplectic cohomology theory  and extend it to  l.c.s. manifolds. This is possible due to  the validity of the Lefschetz
decomposition for these manifolds. The main tool is the spectral sequence
developed by  Di Pietro and Vinogradov for symplectic manifolds, which has
been now adapted and developed further for l.c.s. manifolds. We obtain new results and  applications
even for symplectic manifolds. In
particular we unify various isolated results of the linear symplectic
cohomology theory.

The structure of this  note is as follows.  In section 2 we introduce important linear operators on
l.c.s. manifolds and  study their properties.  The Lefschetz filtration  on the space  $\Om ^* (M^{2n})$  of a l.c.s. manifold $(M^{2n}, \om, \theta)$  is discussed in section 3 together with differential operators  respecting this  filtration.  Then we  use this filtration to construct primitive
cohomology groups  for $(M^{2n}, \om, \theta)$  (Definition \ref{newcomg}).  Some simple properties  of these  groups  are fixed in Proposition \ref{plus1}, Proposition \ref{primi} and Proposition \ref{dpm1}, and their relations with previously  proposed constructions are discussed (Remark \ref{his1}).
The spectral sequence associated with the Lefschetz filtration   is studied in section 4. In particular,  its   $E_1$-term  is compared with the primitive  (co)homological groups (Lemma \ref{fk1}) and the conformal  invariance of this term is proven (Theorem \ref{conf}).  In section 5 we  find some cohomological conditions on $(M^{2n},  \om, \theta)$
under which  this spectral  sequence  stabilizes at  the $E_{t}$-term  (Theorems \ref{stab1}, \ref{nil}, \ref{stcom}).  The last of these theorems gives an answer to the Tseng-Yau question on relations between
 the primitive cohomology and the de Rham cohomology of a compact symplectic manifold.
In section 6 we specialize  the previous theory to  K\"ahler manifolds and prove that for  K\"ahler 
manifolds the spectral sequence  stabilizes  already at its first term (Theorem \ref{kaehler1}).  In section 7  we   compute the primitive  cohomology groups of a compact  $(2n+2)$-dimensional l.c.s. nilmanifold and  a compact 4-dimensional l.c.s. solvmanifold (Propositions \ref{mil1}, \ref{solv1}).  We   study  some properties  of primitive cohomology groups of l.c.s.  manifolds associated with   a co-orientation preserving contactomorphism (Proposition \ref{inv2}). In particular, we show that  the compact l.c.s. solvmanifold is associated with a non-trivial   co-orientation preserving contactomorphism (Theorem \ref{new1}).

The cohomological theory  developed in this note  and its analogues have a much wider area of applications. For instance, it may be naturally adopted to the  class of Poisson symplectic stratified spaces  introduced in \cite{LSV2010}, since these singular  symplectic spaces also enjoy the Lefschetz decomposition.

This   project  was started   as a joint work  of us with Alexandre Vinogradov based on  H.V.L. preliminary  results on l.c.s. manifolds.   Alexandre  Vinogradov   has suggested
 us to extend  the results    to a slightly larger category of twisted  symplectic manifolds. He made  considerable  contributions to improve  the original text  written by H.V.L.,  which we  appreciate   very much.       Eventually  we have noticed that   our  viewpoints   are so different, so we decide  to write  the subject    separately:   in this paper  we deal   only with l.c.s. manifolds and     Alexandre Vinogradov  will deal with  the extension  to  twisted symplectic manifolds.

\section{Basic operators on  a l.c.s. manifold}

In this section we introduce and study basic linear  differential operators  acting on differential forms on a 
l.c.s. manifold $(M^{2n}, \om, \theta)$. 

The first operator  we need is the Lichnerowicz  deformed differential  $d_\theta: \Om^* (M^{2n}) \to \Om^*(M^{2n})$, 
\begin{equation}
 d_\theta (\alpha): =d\alpha + \theta \wedge\alpha.
 \label{lic1}
\end{equation}

Clearly $d_\theta ^2 = 0$ and $ d_\theta (\om) = 0$. The resulting Lichnerowicz cohomology  groups (also called the  Novikov  cohomology groups) are important  conformal invariants  of l.c.s. manifolds. 

Recall that   two  l.c.s. forms $\om$ and $\om'$ on $M^{2n}$ are {\it conformally equivalent}, if
$\om' = \pm (e^ f) \om$ for  some $f \in C^\infty (M^{2n})$.  In this case the corresponding Lee forms $\theta$ and $\theta '$ are cohomologous: $\theta' = \theta \mp df$, hence
 $ d_{\theta}$ and $d_{\theta'}$ are {\it  gauge equivalent}: 
$$d_{\theta'} ( \alpha) = (d_{\theta} \mp df \wedge )\alpha =  e  ^{\pm f} d (e ^{\mp f} \alpha).$$
It follows that  $H^* (\Om^* (M^{2n}), d_\theta)$ and
$H^*(\Om^*(M^{2n}), d_{\theta'})$ are isomorphic. The isomorphism $ I_f : H^* (\Om^* (M^{2n}), d_{\theta}) \to H^*(\Om^*(M^{2n}), d_{\theta'})$ is given by the conformal  transformation $[\alpha] \mapsto [\pm e^ f \alpha]$.

 Note that $d_\theta$ does not satisfy the Leibniz property, unless $\theta = 0$, since
\begin{equation}
d_{\theta} ( \alpha \wedge \beta)  =  d_{\theta} \alpha \wedge \beta + (-1) ^{deg\, \alpha}\alpha \wedge d\beta =  d\alpha \wedge \beta  +(-1) ^{deg\, \alpha}\alpha \wedge d_\theta  \beta. 
\label{lic2n}
\end{equation}

Thus the cohomology group $H^*(\Om^* (M^{2n}), d_{\theta})$  does not have a ring structure, unless $\theta = 0$.
The formula (\ref{lic2n}) also implies that $H^*(\Om^* (M^{2n}), d_{\theta})$  is a $H^*(M, \R)$-module. 

Now let us consider the next basic linear operator
\begin{equation}
L: \Om^* (M^{2n}) \to \Om^* (M^{2n}) , \:  \alpha \mapsto  \om \wedge \alpha. \label{defl}
\end{equation}

Substituting $\alpha := \om$ in (\ref{lic2n})  we obtain   a  nice relation between $d$, $L$ and $d_{\theta}$
\begin{equation}
d_\theta L =  Ld.\label{lic3}
\end{equation}

The identity (\ref{lic3})  suggests us   to consider  a family  of  operators $d_{k\theta}$, which we abbreviate as
$d_k$ if no misunderstanding occurs.  We derive immediately from (\ref{lic3})
\begin{equation}
d _k L^ p = L^p d _{k -p}. \label{lic4}
\end{equation}

The following Lemma is a generalization of (\ref{lic2n}) and it plays an important role in our study of the spectral sequences introduced in later sections. It is obtained by straightforward calculations, so we omit its proof.

\begin{lemma}\label{comd1}  For any $\alpha, \beta \in \Om^* (M^{2n})$  we have
\begin{equation}
d_{k+l} (\alpha \wedge \beta) = d_k \alpha \wedge \beta + (-1) ^{ deg\, \alpha}\alpha \wedge d_l \beta.\label{comd2}
\end{equation}
Consequently
\begin{equation}
d_k \alpha \wedge d_l \beta = d_{k+l} ( \alpha \wedge d_l \beta).\label{comd3}
\end{equation}
Formula (\ref{comd2}) yields  the induced map $H^*(\Om^*(M^{2n}), d_k)  \times H^* (\Om^* (M^{2n}), d_l) \to H^* (\Om^*(M^{2n}), d_{k+l})$.
\end{lemma}

Denote by  $G_{\om}$ the section  of the bundle $\Lambda ^2 TM^{2n}$  such that for all $x \in M^{2n}$ the linear map $i_{G_\om (x)} : T^*_x M^{2n}  \to T_xM^{2n}, \, V  \mapsto i_V (G_\om (x)),$  is the inverse of the map $I_\om : T_x M^{2n}  \to T^*_xM^{2n},\, V \mapsto i_V \om$. 
Clearly $G_\om$  defines a bilinear pairing: $T^*M^{2n} \times T^*M^{2n}\to C^\infty (M^{2n})$. Denote by $\Lambda ^p G_\om$ the associated pairing: $\Lambda^p (T^*M) \times \Lambda ^p(T^*M) \to C^\infty (M^{2n})$. The   l.c.s.  form $\om$ and the   associated bi-vector field $G_\om$  define {\it a l.c.s. star  operator} $*_\om: \Om^p (M^{2n}) \to \Om ^{2n-p}(M^{2n})$  as follows  \cite[\S 2.1]{Brylinski1988}.
\begin{equation}
 *_\om : \Om ^p (M^{2n}) \to \Om^{2n-p} (M^{2n}),\,  \beta \wedge *_\om \alpha : = \Lambda ^p G_\om  (\beta, \alpha) \wedge \frac {\om ^n} { n !},
\label{sst}
\end{equation}
for all $\alpha, \beta  \in  \Om ^p (M^{2n})$.
Using \cite[Lemma 2.1.2]{Brylinski1988} we get easily
\begin{equation}
*_\om ^2  = Id.\label{bryss1}
\end{equation}
We  define  {\it the l.c.s. adjoint} $L^*$  of $L$ and { \it  the l.c.s. adjoint} $(d_k)^*_\om$  of $d_k$ with respect to  the  l.c.s. form $\om$ as follows  
\begin{equation}
L^* : \Om ^{p} (M^{2n})\to \Om^{p-2} (M^{2n}), \, \alpha^p \mapsto - *_\om L *_\om \alpha^p. \label{defls}
\end{equation}
\begin{equation}
(d_k)^*_\om: \Om ^{p} (M^{2n})\to \Om^{p-1} (M^{2n}), \,  \alpha^p \mapsto (-1) ^p * _\om d_{n+k-p} *_\om (\alpha^p). \label{defnbls}
\end{equation}
For symplectic manifolds our definition of $(d_k)^*_\om$ agrees with the one in  \cite[\S 1]{Yan1996},  it is different from the one in \cite[Theorem 2.2.1]{Brylinski1988} by 
 sign (-1). 

A  section $g $ of the bundle $S^2 T^*M^{2n}$ is called {\it a compatible   metric}, if  there is an almost complex structure $J$ on $M^{2n}$ such  that $g (X, Y) = \om(X, JY)$. In this case $J$ is called {\it a compatible almost complex structure.}
Recall that  the Hodge operator  $*_g$ is defined as follows
\begin{equation} 
*_g : \Om ^p (M^{2n}) \to \Om^{2n-p} (M^{2n}),\,  \beta \wedge *_g \alpha : = \Lambda ^p G_g  (\beta, \alpha) \wedge \frac {\om ^n} { n !},\label{hodgets}
\end{equation}
where $G_g \in \Gamma (S^2 TM^{2n})$ is the ``inverse of $g$", i.e. it is defined in the same way as  we define $G_\om$ above: for all $x \in M^{2n}$ the linear map $i_{G_g (x)} : T^*_x M^{2n}  \to T_xM, \, V  \mapsto i_V (G_g (x)),$  is the inverse of the map $I_g : T_x M  \to T^*_xM,\, V \mapsto i_V (g)$. We also denote by $\Lambda  ^p G_g$   the associated pairing: $\Lambda^p (T^*M)  \times \Lambda ^p(T^*M) \to C^\infty (M^{2n})$  induced by $G_g$, (see also \cite[p.105]{Brylinski1988} for comparing $\Lambda ^p G_\om$  with $\Lambda^p G_g$).

Using \cite[Lemma 5.5]{Voisin2007} we get easily
\begin{equation}
*_g ^ 2 (\alpha ^p) =  (-1) ^p\alpha ^p \text{ for } \alpha ^p \in \Om ^p (M^{2n}). \label{voisin5.5}
\end{equation}

\begin{lemma}\label{fund} 1. The space of  metrics  compatible with a given l.c.s.  form  $\om \in \Om^2 (M^{2n})$ is contractible.

2. (cf.  \cite[chapter II, 6.2.1]{Voisin2007}) In the presence of a compatible  metric $g$ on $M^{2n}$  we have
\begin{equation}
L^* =  \Lambda \label{lel}
\end{equation}
where $\Lambda = (*_g)^{-1} L *_g$ is the  adjoint of $L$ with respect to the metric $g$.
\end{lemma}

\begin{proof} 1. The proof  for the first assertion goes in the same way  as for the case  of symplectic manifolds, so we omit its proof.

2. The second assertion of Lemma \ref{fund} is a simple consequence of  the following 

\begin{lemma}\label{hk}\cite[Theorem 2.4]{Brylinski1988} Assume that $(M^{2n}, J, g)$ is an almost Hermitian manifold and $\om$ is the associated almost symplectic form.  For $\alpha \in \Om ^{p, q}(M^{2n})$ we have
$$*_\om (\alpha) = \sqrt{-1}^{p-q} *_g (\alpha).$$
Here we extend  $*_\om$ and $*_g$  $\C$-linearly on $\Om ^* (M^{2n})\otimes \C$.
\end{lemma}
This completes the proof of Lemma \ref{fund}.
\end{proof}

Let $\pi_{k} : \Om^*  (M^{2n}) \to \Om^k (M^{2n})$ be  the projection. Denote $\sum _{i=0} ^{2n} (n-k)\pi_{k}$  by $A$. 
Using  well-known identities in K\"ahler  geometry for $(\Lambda, L, A)$, see e.g. \cite[(IV), chapter I]{Weil1958}, \cite[p.121]{GH1978}, \cite[Lemma 6.19]{Voisin2007},  Lemma \ref{fund} implies  immediately the  following
 
\begin{corollary}\label{hks} (cf. \cite[\S 1]{Lychagin1979}, \cite[Corollary 1.6]{Yan1996}) On any  l.c.s. manifold $(M^{2n}, \om, \theta)$ we have
\begin{equation}
L ^*  = i (G_\om),\label{lss}
\end{equation}
\begin{equation}
[L ^*, L] = A,\: [A, L] = - 2L, \: [A, L^*]=  2L^*.\label{sl2}
\end{equation}
\end{corollary}

The relation in (\ref{sl2}) shows that $(L^*, L, A)$ forms a $\ssl_2$-triple, which has many important consequences
for l.c.s. manifolds.

\begin{proposition}\label{comdn} The following commutation relation hold
\begin{equation}
L^* (d _k)^*_{\om}=  (d _{k-1})^*_\om L^*.\label{lic4s}
\end{equation}
\end{proposition}
\begin{proof}  Clearly (\ref{lic4s}) is  obtained from (\ref{lic4}) by applying  the LHS and RHS of (\ref{lic4}) the  l.c.s. star operator  on the left and on the right, taking into account (\ref{bryss1}).
\end{proof}

\section{Primitive forms and  primitive cohomologies}

In this section we introduce the notions of   primitive forms and  coeffective forms on a  l.c.s. manifold
$(M^{2n},  \om, \theta)$, using the linear operators $L$ and $L^*$ defined in the previous section.  As in the symplectic case we obtain   a Lefschetz decomposition  of the space $\Om^*(M^{2n})$  induced  by
 primitive forms and  coeffective forms together with various linear differential operators respecting
this  decomposition as well as an associated filtration (Propositions \ref{clo2} and  \ref{hlc}).
The natural splitting of  the introduced differential operators according to the Lefschetz decomposition leads to new  cohomology groups
of $(M^{2n},\om, \theta)$ (Definition \ref{newcomg}). In  Propositions \ref{plus1},  \ref{primi}, \ref{dpm1} we fix simple properties of these new cohomology groups.
At the end of this section we compare our  construction  with related  constructions in \cite{Lychagin1979}, \cite{Bouche1990}, \cite{Rumin1994},  \cite{FIL1996}, \cite{Yan1996},   \cite{PV2006}, \cite{Pietro2006}.

\begin{definition}\label{primco} (\cite{Bouche1990}, \cite{Yan1996}, cf. \cite{Weil1958}, \cite{GH1978})  An element  $\alpha \in \Lambda ^k T^*_xM^{2n}$, $0\le k \le n$,  is called {\it primitive} (or {\it  effective}), if
$L^{n-k+1} \alpha  = 0$.   An element $\alpha \in \Lambda ^k T^*_xM^{2n}$, $n+1\le k \le 2n$,  is called {\it primitive}, if $\alpha = 0$.  An element $\beta \in \Lambda^k T^*_x M^{2n}$ is called {\it coeffective}, if
$L \beta = 0$.
\end{definition}

\begin{remark}\label{inv1}  1. Wells in  \cite{Wells1986} refers to  Lefschetz \cite{Lefschetz1924}  and    Weil \cite{Weil1958}  for the terminology ``Lefschetz decomposition" and ``primitive forms".  Many 
mathematicians prefer  ``Lepage decomposition" and  ``effective forms"  following Lepage in \cite{Lepage1946}.

2.  Clearly the notion of  primitive form  as well as the notion of coeffective form
depends only on the conformal class of a l.c.s. form $\om$.
\end{remark}

The  relation (\ref{sl2}) between linear operators $L, L^*$ and $A$ leads  to Lemma \ref{yanpr} below characterizing primitive forms and coeffective forms. The  resulting Lefschetz decomposition  of the space $\Lambda  T^* M^{2n}$ 
is a direct consequence of the  $\ssl(2)$-module theory.  Various  variants of Lemma \ref{yanpr} for symplectic manifolds   appeared  in many works, beginning   possibly with  a paper by Lepage \cite{Lepage1946}, with later applications in K\"ahler geometry \cite{Weil1958}, \cite{GH1978}, \cite{Voisin2007}, in a theory  of second-order  differential equations \cite{Lychagin1979},  in symplectic  geometry \cite{Bouche1990}, \cite{Yan1996},  etc..

We  denote by $P^k_x(M^{2n})$ the set of primitive elements in  $\Lambda^k T^*_xM^{2n}$.

\begin{lemma}\label{yanpr}   1. An element  $\alpha \in \Lambda ^k T^*_xM^{2n}$ is primitive, if and only if 
$L^* \alpha = 0$.\\
2. 
 An element $\beta \in \Lambda ^k T^*_xM^{2n}$ is  coeffective, if and only
if $*_\om  \beta $ is  primitive.\\
3. We have the following Lefschetz decomposition for  $n \ge k \ge 0$:
\begin{equation}
\Lambda ^{n-k}T^*_x M^{2n}= P^{n-k }_x (M^{2n}) \oplus LP^{n-k-2}_x(M^{2n}) \oplus  L^2 P^{n-k-4}_x(M^{2n}) \oplus \cdots,
\label{lh1}
\end{equation}
\begin{equation}
\Lambda ^{n+k}T^*_x M^{2n}=L^k P^{n-k }_x (M^{2n}) \oplus L^{k+1}P^{n-k-2}_x(M^{2n})  \oplus \cdots .
\label{lh2}
\end{equation}
\end{lemma}

From Lemma \ref{yanpr}  we get immediately

\begin{corollary}\label{inj1}
1. $L^k : \Lambda^{n-k} T^*_xM^{2n} \to \Lambda ^{n+k} T^*_xM^{2n}$ is an isomorphism, for $0\le k \le n$.\\
2. $L : \Lambda^{n-k-2} T^*_xM^{2n}\to \Lambda^{n-k}T^*_xM^{2n}$ is injective, for $k =-1, 0, 1, \cdots , n-2$.
\end{corollary}

It is useful to introduce the following notations.
Denote by $P^{n-k}M^{2n}$  the subbundle  in $\Lambda T^*M^{2n}$ whose fiber  is $P^{n-k} _x (M^{2n})$.
Let  $\Pp^{n-k}(M^{2n})\subset \Om^{n-k} (M^{2n})$ be the  space of all smooth $(n-k)$-forms with values in  $P^{n-k}M^{2n}$.  
Elements of $\Pp^{n-k}(M^{2n})$
are called   {\it primitive $(n-k)$-forms}.  Let us  set (cf. \cite{TY2009})
\begin{equation}
\Ll ^{s, r} : = L^s \Pp^r  \text { for }  0 \le s, r \le n. \label{lty1}
\end{equation}
Put $\Pp^* (M^{2n}): = \oplus _r \Pp ^r (M^{2n})$. Then Lemma \ref{yanpr} yields the following  decompositions, which we call  {\it the first and  second Lefschetz decompositions}
\begin{equation}
\Om^*(M^{2n}) =\Pp^*(M^{2n}) \oplus L \Om^* (M^{2n}) =  \bigoplus _{0\le 2s + r \le 2n} \Ll ^{s, r}.\label{lty2}
\end{equation}

Now we  consider the  interplay  between the Lefschetz decompositions
(\ref{lty2}) and  the linear differential operators introduced in the previous section. 
Iterating the action of $L$ on $K^*: = \Om ^* (M^{2n})$, we define the following filtration 
\begin{equation}
F^0 K^*: = K ^*  \supset  F^1 K^*: =   LK^* \supset   \cdots 
\supset F^k K^* :=L^{k} K^*\supset \cdots    \supset F^{n+1} K^* =
 \{ 0\} .
\label{vin2}
\end{equation}

\begin{proposition}\label{clo2} 1. The subset $F^k K^* $ is stable  with respect to  $d_p$  for all $k$ and $p$.\\
2. For any $\gamma \in \Om ^1 (M^{2n})$ we have 
$$\gamma \wedge \Ll^{0, n-k}  \subset \Ll^{0, n-k+1}\oplus \Ll^{1,n-k-1}.$$
\end{proposition}

\begin{proof} 1.  The first assertion of Proposition \ref{clo2} follows from the identity $d_p (\om ^k \wedge \phi) = \om^k\wedge d_{p-k} \phi$ for $\phi \in \Om^*(M^{2n})$.\\
2.  Assume that $\alpha \in \Pp^{n-k}(M^{2n}) = \Ll^{0,n-k}$. Then  $L^{k +1} (\gamma \wedge \alpha) = \gamma \wedge L^{k+1} (\alpha) = 0$.  Taking into account  the decomposition of $\gamma\wedge \alpha$  according to the second Lefschetz decomposition we obtain  the second assertion of  Proposition \ref{clo2} immediately.
\end{proof}

\begin{remark}\label{rem:left1} 1. The relation $d_p(\om ^k \wedge \phi)  = \om ^ k \wedge d_{p-k} \phi$   can be  also interpreted  as an interplay between  different filtered complexes $(F^*K^*, d_k)$ and
$(F^*K^*, d_p)$. We shall investigate this interplay  deeper   in the next section.

2. We observe that the  decompositions (\ref{lh1}), (\ref{lh2})  and (\ref{lty2})  are compatible with the filtration (\ref{vin2}) in the following sense.  For any  $p\ge 0$  and $0 \le k \le 2n$  we have  
\begin{equation}
F^p K^* \cap  \Om ^{k} (M^{2n})  = \oplus_{i =0}^{[{k\over 2}]-p} \Ll^{p+i, k -2p-2i} \text{ if } k \ge 2p,
 \label{fil1}
\end{equation}
\begin{equation}
 F^p K^*   \cap  \Om ^{k} (M^{2n})  = 0
 \text{ if }   k < 2p,
 \label{fil1a}
\end{equation}
The decomposition in  (\ref{fil1}) and (\ref{fil1a})   will be called {\it the induced  Lefschetz decomposition}. It is important for understanding  the  spectral sequences  introduced in the next section.
\end{remark}

\begin{proposition}\label{hlc} The following inclusions hold  
\begin{equation}
d_r\Ll^{p, q-p} \subset  \Ll^{p,q-p+1}\oplus  \Ll^{p+1,q-p-1},
\label{prdlt}
\end{equation}
\begin{equation}
(d_r)^*_\om  \Pp^{n-k}(M^{2n})\subset \Pp^{n-k-1} (M^{2n}).
\label{prde1}
\end{equation}
\end{proposition}

\begin{proof}   Let $\beta\in \Pp^{q} (M^{2n})= \Ll ^{0,q}$, so $L^{n-q+1}\beta  = 0$.  We derive  from (\ref{lic4}) 
\begin{equation}  
L^{n-q+1} d _r\beta  = d_{r+n - q +1}L^{n-q+1} \beta =  0.
\end{equation}
 Using  (\ref{lh1}) and (\ref{lh2})
we get $d_{r} \beta \in \Pp^{q+1}(M^{2n}) + L\Pp^{q-1}(M^{2n})$. This  proves the  inclusion (\ref{prdlt})  of Proposition \ref{hlc}  for $p = 0$.
The  inclusion (\ref{prdlt}) for $p \not = 0$ follows from the  particular case  $p = 0$ and the identity $d_r L^p = L^p d_{r-p}$.

 Assume that $\beta \in \Pp^{n-k} (M^{2n})$. Taking into account   (\ref{lic4s})  we  obtain 
$$ L^*(d_r)^*_\om \beta  =(d _{r-1})^*_\om L ^*  \beta,$$
which is zero since $\beta$ is primitive.
Hence $(d_r)^*_\om \beta$ is also primitive. This proves  (\ref{prde1}) and completes the proof of Proposition \ref{hlc}.
\end{proof}

Now we will show several consequences of Proposition \ref{hlc}. Denote by $\Pi_{pr}$ the projection $\Om^* (M^{2n}) \to \Pp^*(M^{2n})$ according to the Lefschetz decomposition in (\ref{lty2}).
Set 
$$d_k^+:= \Pi_{pr}d_k.$$
Using the first Lefschetz decomposition   and Proposition \ref{hlc}    we decompose the operator  $d_k: \Om ^q (M^{2n}) \to \Om ^{q +1} (M^{2n})$ for $ 0\le q \le n$ as follows (cf. \cite{TY2009}). 
\begin{equation}
d_k = d_{k}^+ + L d_{k}^-,\label{hoya}
\end{equation}
where  $d^-_{k} : \Om ^q (M^{2n}) \to \Om ^{q-1} (M^{2n}), \, 0\le q \le n$.  Note that $ d ^-_k$ is well-defined, since $L : \Om ^{q-1} (M^{2n}) \to \Om ^{ q+1}(M^{2n})$ is
injective.  It is straightforward to check
\begin{equation}
d_k^+(\Ll ^{s, r}) = 0 \text { if }  s\ge 1, \text{ and } d_{k}^-(\Ll ^{s, r}) \subset \Ll ^{s, r-1}.\label{lty3}
\end{equation}

\begin{lemma}\label{hoya1} (cf. \cite[Lemma 2.5, II]{TY2009}) The operators $d_k^+$, $d_{k-1}^-$   satisfy the following properties
\begin{eqnarray} (d_k^+) ^2 (\alpha^q) = 0,\label{com0}\\
d_{k-1} ^-d_k^-(\alpha ^q)= 0 , \text { if } q \le n,  \label{com1}\\
(d_k^- d_k^+ + d_{k-1}^+ d^-_k) \alpha ^ q = 0, \text{ if } q \le n-1, \label{com1a}\\
(d_{k-1}) ^* _\om (d_k) ^* _\om (\alpha^q) = 0.\label{prag0510}
\end{eqnarray}
\end{lemma}
\begin{proof} We use the equality $d_k^2=0$ in the form $d_k(d_k^++Ld_k^-)=0$. Using (\ref{lic4}) we get
$$
(d_k^+)^2+L(d_k^-d_k^++d_{k-1}^+d_k^-)+L^2d_{k-1}^-d_k^-=0.
$$ 
Now taking into account (\ref{lty3})  and the injectivity of the operators $L : \Om ^{q} (M^{2n}) \to  \Om ^{q+2}(M^{2n})$ and $L ^2 : \Om ^{q-1} (M^{2n}) \to \Om ^{q+3} (M^{2n})$ for $q \le n-1$ we obtain (\ref{com0}), (\ref{com1}), and (\ref{com1a}). 

Finally, (\ref{prag0510}) is a consequence of $d_k^2=0$ and $*_\omega^2=Id$. 
\end{proof}

\medskip

Proposition  \ref{hlc} and Lemma \ref{hoya1} lead to  new cohomology groups   associated with a l.c.s. manifold $(M^{2n}, \om, \theta)$.  We observe that $\Pp^*(M^{2n})$ is stable under the action of the operators $d_k^+,  (d _k)^*_\om, d_k ^-$. 

\begin{definition}\label{newcomg}  Assume that $0\le q \le n-1$.\\
The {\it $k$-plus-primitive  q-th cohomology group}  of $(M^{2n}, \om, \theta)$ is defined by
\begin{equation}
H^{q}  (\Pp^*(M^{2n}), d_k ^+): = \frac{\ker d_{k} ^+ :\Pp^{q} (M^{2n}) \to \Pp ^{q+1}(M^{2n})}{ d_{k}^+ (\Pp^{q-1} (M^{2n}))}.\label{vy3}
\end{equation}
The {\it $k$-primitive q-th-cohomology group} of $(M^{2n}, \om, \theta)$ is defined by
\begin{equation}
H^{q}   (\Pp^*(M^{2n}), (d_k)^*_\om) : = \frac{\ker (d_k)_\om^*: \Pp^{q} (M^{2n})\to \Pp ^{q-1}(M^{2n})}{ (d_{k+1})_\om^* (\Pp^{q+1}(M^{2n}))}.\label{vy3n}
\end{equation}
 The {\it $k$-minus-primitive q-th cohomology group} of $(M^{2n}, \om, \theta)$ is defined by
\begin{equation}
H^{q}(\Pp^*(M^{2n}), d _k ^-):= \frac{\ker d_k^- : \Pp^q(M^{2n})\to \Pp^{q-1}(M^{2n})}{d_{k+1}^- ( \Pp^{q+1}(M^{2n}))}.\label{vy4}
\end{equation}
\end{definition}

Now we show few simple properties of  the associated  cohomology groups  of  a l.c.s. manifolds.  
Note that the formula (\ref{plus1b}) below   has been proved in
\cite[Proposition 3.15]{TY2009} for  compact symplectic manifolds $(M^{2n},\om)$.

\begin{proposition}\label{plus1}   Assume that $(M^{2n}, \om, \theta)$ is a l.c.s. manifold, $n \ge 2$. 

1. Suppose  that $[(k-1)\theta]\not = 0\in H^1 (M^{2n}, \R)$. Then  
\begin{equation}
H^{1}(\Pp^*(M^{2n}), d_k ^+) =  H^1  (\Om^*(M^{2n}), d_{k}). \label{plus1a}
\end{equation}

2.  Suppose that  $[(k-1)\theta] = 0\in H^1 (M^{2n}, \R)$.  Then
\begin{equation}
H^{1} (\Pp^*(M^{2n}), d _k ^+) = H^1(\Om^*(M^{2n}), d_\theta) \text { if } [\om] \not = 0 \in H^2(\Om^*(M^{2n}), d_\theta),\label{plus1b}
\end{equation}
\begin{equation}
H^{1} (\Pp^*(M^{2n}), d_k ^+) =  H^1(\Om^*(M^{2n}), d_\theta)\oplus R \text { if } [\om]  = 0 \in H^2 (\Om^*(M^{2n}), d_\theta),\label{plus1c}
\end{equation}
where $R$ is the 1-dimensional vector space generated by  $\rho \in H^{1} (\Pp ^*(M^{2n}), d_k ^+)$ with $d _{k} \rho = \om$. 
\end{proposition}

\begin{proof} 1. Assume that  $0\not=\alpha \in \Pp^1(M^{2n})$ and $d_k ^+ \alpha = 0$, i.e. $[\alpha] \in  H^{1} (\Pp ^*(M^{2n}), d_k ^+)$. Since $d^+_k \alpha = 0$ we get
$d_k \alpha = L f$,  where $f\in C^\infty (M^{2n})$.   Assume that $f \not = 0$.
Using $d^2_k \alpha = 0$ we derive $Ld_{k-1}f = 0$, which implies $d_{k-1}f =0 $, since $L$ is injective. The equality $d_{k-1} f = 0$ implies   that $d_{k-1}$  is  gauge  equivalent to $d$. This contradicts 
the assumption of Proposition \ref{plus1}.1. Hence  $f = 0$. It follows  $d _k \alpha = 0$. Using $d ^+_kh = d _k h$ for all $h \in \Pp^0 (M^{2n})$ we
obtain  (\ref{plus1a}) immediately.

2. Now we assume that $[(k-1)\theta] = 0\in H^1 (M^{2n}, \R)$, or equivalently, $d_{k-1}$ is gauge equivalent to the canonical connection : $d_{k-1} =  e^h d e^{-h} = d - dh\wedge$  for some $h \in C^\infty (M ^{2n})$. In this case, as above, $d _{k}\alpha  = Lf$ implies   $d_{k-1}f = 0$, and hence
$f = ce^{h}$.  

a)  Assume  that (\ref{plus1b}) holds. If $c \not = 0$, then     $[e^h\om] = 0 \in  H^2(M^{2n}, d_k)$.  Since $(k-1)\theta = -dh$, the deformed differential  $ d _\theta - dh\wedge$  is   gauge equivalent to $d {_k}$.
It  follows that   $[\om] = 0 \in H^2 (M^{2n}, d_\theta)$. This contradicts the assumption  of (\ref{plus1b}).  
 Hence  $c =0$.  In this case  we have $d_{k} \alpha = 0$, and therefore $[e^{-h}\alpha] \in H^1(M^{2n}, d_\theta)$. Taking into account $d_k^+ f = d_k  f$ for any
 $f \in \Pp^0 (M^{2n}) = \Om^0 (M^{2n})$,  we obtain  (\ref{plus1b}).
 
b)  Assume that (\ref{plus1c}) holds.  Then $e^h\om = d_k \rho$ for some $\rho \in \Om^1 (M^{2n})$.  
In this case, $d_k(\alpha)  =  L f = \om  c e ^h $ implies  $d_k (\alpha - c\rho) = 0$. We conclude that  if $d_k ^+(\alpha) = 0$ then $\alpha = c \rho + \beta$ where
$d_k (\beta) =0$. Clearly $[\beta] \in H^1(M^{2n}, d_k) = H^1(M^{2n}, d_\theta)$.
This proves (\ref{plus1c}) and completes the proof of Proposition \ref{plus1}.
\end{proof}

\begin{proposition}\label{primi} (cf. \cite[Lemma 2.7, part II]{TY2009}) Assume that  $0\le k \le n$.
 If  $\alpha \in \Pp^k (M^{2n})$, then 
\begin{equation}
d^-_r (\alpha ^k) = {(d_r)^*_\om (\alpha ^k)\over  n-k+1}. \label{plus13}
\end{equation}
Consequently $H^{k} (\Pp^*(M^{2n}), d^-_r) = H^{k}(\Pp^*(M^{2n}), (d_r) ^* _\om)$.
\end{proposition}

\begin{proof} It suffices to prove (\ref{plus13}) locally.  Note that locally $d _\theta = d - df\wedge$. In this case $(d-df\wedge)\om = 0$ implies $\om = e^f \om_0$ with $d \om_0 = 0$. Next, we compare $*_\om$ and $*_{\om_0}$, using (\ref{sst}) and the equality $G_\om = e ^{-f}G_{\om_0}$.
$$\beta^k \wedge *_\om \alpha^k  = \wedge ^k G_\om  (\beta^k, \alpha^k) \wedge \frac {\om ^n} { n !} = \wedge ^k e ^{-kf}G_{\om_0}(\beta^k, \alpha  ^k) e^{nf}\frac {\om_0 ^n} { n !}= e^{(n-k)f} \beta^k\wedge*_{\om_0} \alpha^k,$$
where $\beta^k , \alpha ^k \in  \Om ^k (M^{2n})$. It follows that
\begin{equation}
* _\om(\alpha ^k)  = e ^{(n-k)f}*_{\om_0} (\alpha ^k).\label{laho2}
\end{equation}

Let $\alpha^k \in \Pp^k (M^{2n})$, $0\le k \le n$. Denote by $(d)^*_{\om_0}$ the symplectic adjoint of $d$ with respect to $\om_0$. The formula (\ref{plus15}) below, which is a partial case of (\ref{plus13}) for  symplectic manifold, has been proved in \cite[Lemma 2.7, part II]{TY2009}. (We observe that their   operator  $d^\Lambda$ differs from our  operator $(d)^*_{\om_0}$ by sign (-1).)
\begin{equation}
d ^- (\alpha ^k) = { (d )^*_{\om_0} \alpha ^k \over n-k+1}.\label{plus15}
\end{equation}
 Using $d_{n+r-k}\alpha = e^{(n+r-k)f}d(e^{-(n +r-k)f}\alpha)$ we obtain from (\ref{laho2})
$$
(d_r)^*_\om (\alpha^k)= (-1) ^k*_\om d_{n+r-k}*_\om (\alpha ^k) =$$
$$= (-1)^k e^{(n - (2n -k+1))f}*_{\om_0} e^{(n+r-k)f}d(e^{-(n +r-k)f}  ( e^{(n-k) f}*_{\om_0} \alpha ^k)=$$ 
\begin{equation}
= (-1) ^k e ^{(r-1) f}*_{\om_0}d ( e^{ -rf}*_{\om_0} \alpha ^k) =\label{gauge1}
\end{equation}
\begin{equation}
 = e ^{-f} ( (d)_{\om_0}^* \alpha ^k + (-1) ^k *_{\om_0}(-r) df \wedge *_{\om_0}\alpha ^k). \label{lahor6}
\end{equation}
Substituting $r = 0$,  we derive  from (\ref{lahor6})
\begin{equation}
(d)_{\om_0} ^* (\alpha ^k) = e^f (d) ^*_\om\alpha ^k . \label{laho5}
\end{equation}
Next we  compare $d_r ^-$  with $d^-$.
\begin{equation}
d _r (\alpha^k) = e^{rf} d (e ^{-rf} \alpha^k) = e^{rf} [d^+ (e ^{-rf} \alpha ^k) + \om_0 \wedge d^- e^{-rf} \alpha ^k].\label{laho3}
\end{equation}
 Since  the Lefschetz decomposition of $\Om(M^{2n})$ is invariant under  conformal transformations we obtain from (\ref{laho3})
 \begin{equation}
 d_r^-  (\alpha ^k) = e^{(r-1)f} d^- (e ^{-rf} \alpha ^k).\label{laho1}
 \end{equation}
 Combining (\ref{laho1}) with (\ref{plus15}) we conclude that
\begin{equation}
d_r ^- (\alpha ^k)=  e^{(r-1)f}\frac{(d)_{\om_0}^* (e^{-rf} \alpha ^k)}{n-k+1}.\label{laho7}
\end{equation}
Taking into account  (\ref{laho5}) and (\ref{gauge1}), we derive from (\ref{laho7})
\begin{equation}
d^-_r (\alpha ^k) = e^{rf}\frac{(d)^*_{\om}e^{-rf} \alpha ^k}{ n-k+1} = \frac{(d_r)^*_\om \alpha ^k}{n-k+1} .\label{plus16}
\end{equation}
This proves (\ref{plus13}). 

Clearly the second assertion of Proposition \ref{primi} follows from (\ref{plus13}). This completes the proof of Proposition \ref{primi}.
\end{proof}

Let $J$ be a compatible almost complex structure on a  l.c.s. manifold $(M^{2n}, \om, \theta)$. The complexified space $T^*_\C (M^{2n}) : =(T ^*(M ^{2n})\otimes \C $  is decomposed into
eigen-subspaces $T^{p,q}(M^{2n})$. Let $\Pi^{p,q}: T_\C^*(M^{2n}) \to  T^{p,q} (M^{2n})$ be the projection.
Set
$$\Jj : = \sum _{p,q} (\sqrt{-1})^{p-q} \Pi ^{p,q}.$$

In what follows we want to apply the Hodge theory for compact l.c.s. manifolds $(M^{2n}, \om, \theta)$ provided with a compatible   metric $g$.  First we derive a formula for the formal  adjoint $(d ^+ _l)^*$ of $d ^+  _l: \Pp^*(M^{2n}): = \Pp ^* (M^{2n}) \to \Pp^*(M^{2n})\subset \Om^*(M^{2n})$.   For any operator $D$ acting on a subbundle $E \subset \Om^* (M^{2n})$ we denote by $(D)^*$ the  formal adjoint of $D$.

\begin{lemma}\label{formal1} For any $\alpha \in \Pp^* (M^{2n})$ we  have     
\begin{eqnarray}
(d^+ _l)  ^*  (\alpha) = - *_g (d  _{-l}) *_g(\alpha). \label{formal3}
\end{eqnarray}
\end{lemma}

\begin{proof}  First,  we want to compute the  formal adjoint $(d_l) ^* $  of $d _l = d + l \theta \wedge : \Om ^* (M^{2n}) \to \Om ^* (M^{2n})$. It is known that  \cite[\S 5.1.2]{Voisin2007}
\begin{equation}
(d)^* = - * _g  d *_g.\label{voisin521}
\end{equation} 

Since $\theta \wedge$ is the symbol of $d$ we derive  from (\ref{voisin521})
\begin{equation}
(l\theta \wedge ) ^*  = *_g l\theta \wedge *_g. \label{formal7}
\end{equation}
It follows from (\ref{voisin521}) and (\ref{formal7})
\begin{equation}
(d _l )^*= -*_gd _{-l}*_g.\label{formal4}
\end{equation} 
Using (\ref{lel})  we get
for $\alpha \in \Pp ^* (M^{2n})$ 
\begin{equation}
(Ld_l ^-) ^*( \alpha) = (d_l ^-)^* \Lambda (\alpha) = 0.\label{formal9}
\end{equation}
It follows from (\ref{formal4}) and (\ref{formal9}) that for $\alpha \in \Pp ^* (M^{2n})$
\begin{equation}
(d_l^+) ^*(\alpha) =  -*_g (d_{-l}) *_g(\alpha) . \label{formal10}
\end{equation}
This proves (\ref{formal3}), which   is consistent  with \cite[(3.2), part II]{TY2009}, if $(M^{2n}, \om)$ is a symplectic manifold.
\end{proof}

\begin{lemma}\label{dpm} (cf. \cite[Lemma 3.4, part II]{TY2009}) Let $J$ be a compatible almost complex structure on a   l.c.s.  manifold  $(M^{2n}, \om, \theta)$, $g$  the  associated compatible   metric and $*_g$ the   Hodge  star operator with respect to $g$.
Then for  $\alpha ^k, \alpha ^{k-1} \in \Pp^* (M^{2n})$, $0\le k \le n$, we have
\begin{eqnarray}
\Jj  (d^+_l) ^*  \Jj^{-1}(\alpha^k)  =  (n-k+1) d^- _{-l+k -n}(\alpha^k),\label{conj2}\\
\Jj  d ^+_l \Jj^{-1} (\alpha ^{k-1}) = (n-k+1) (d^- _ {-l+k-n}) ^*(\alpha^{k-1}). \label{conj1}
\end{eqnarray}
\end{lemma}

\begin{proof} Using \cite[Theorem 2.4]{Brylinski1988}, see also Lemma \ref{hk},  we get easily 
\begin{equation}
\Jj  = *_g *_\om. \label{conj3}
\end{equation}
 By (\ref{voisin5.5}) we derive from (\ref{conj3})
\begin{equation}
\Jj ^{-1} (\alpha ^k)  = *_\om *_g (-1) ^k(\alpha^k).\label{conj3s}
\end{equation}
 Combining (\ref{conj3s}) with (\ref{conj3}),  (\ref{formal3}) and applying (\ref{defnbls}), (\ref{voisin5.5}) again we obtain
\begin{eqnarray}
\Jj (d^+ _l)^*_g \Jj ^{-1} (\alpha ^k) = (-1) ^{k+1} *_g *_\om *_g d _{-l} *_g *_\om *_g(\alpha^k)=\nonumber\\
 = (-1)^{k+1}*_g ^2   ( d_{-l+k -n}) ^*_\om (\alpha^k) =  (d_{-l +k -n})^*_\om (\alpha^k),\label{conj4}
\end{eqnarray}
since $*_\om *_g   = *_g *_\om$.
Using (\ref{plus13}) we derive (\ref{conj1}) immediately from  (\ref{conj4}).  Clearly (\ref{conj1})  follows from (\ref{conj2}), since they are  adjoint. This completes the proof of  Lemma \ref{dpm}.
\end{proof}

The following Proposition is a generalization of  \cite[Proposition 3.5, part II]{TY2009}.

\begin{proposition}\label{dpm1} Let  $(M^{2n},  \om, \theta)$ be a compact l.c.s manifold. Then there is  $H ^{k} (\Pp ^*(M^{2n}), d^+_l) =  H^{k} (\Pp ^*(M^{2n}),(d_{-l+k -n})^*_\om )$ for all  $ l \in \Z$ and $0\le k \le n-1$.
\end{proposition}

\begin{proof} First we note that all the operators $d ^+_l$, $d^-_l$ and $(d_l)^* _\om$ restricted to the  space $\Pp^*(M^{2n})$ are elliptic operators.  This  observation is a consequence of 
the  theorem by Bouche who proved that the  complex of coeffective forms on a symplectic manifold  $M^{2n}$ is elliptic in dimension greater than $n$ \cite{Bouche1990}. Indeed,  the complex $(\Pp^*(M^{2n}), (d)_\om ^*)$ is dual to the complex
of  coeffective forms, see also Remark \ref{his1}.1 below.  Thus  $(d_l) ^*_\om$  acting on $\Pp ^*(M^{2n})$ is   an elliptic operator, since  $(d_l)^*_\om$ has the same symbol as $(d)_\om^*$.  Taking  (\ref{plus13}), (\ref{conj1}) and (\ref{conj2}) into account  we prove the ellipticity of $d^-_l$ and $d^+_l$ acting
on $\Pp^*(M^{2n})$.  In \cite[Proposition 2.8 part II]{TY2009} the authors give another proof of the ellipticity of these operators.

Now Proposition \ref{dpm1} follows easily  from Lemma \ref{dpm} using the Hodge theory. 
\end{proof}

\begin{corollary}\label{lich0}  Assume  that $(M^{2n}, \om, \theta )$ is a connected compact  l.c.s. manifold. Then $H^{0} (\Pp^*(M^{2n}), d^+_k) = 0 $ if $d_k$  is not gauge equivalent to the
canonical differential $d = d_0$, ore equivalently $[k\theta] \not = 0\in H ^1 (M^{2n}, \R)$.
If otherwise, then $H^{0} (\Pp^*(M^{2n}), d^+_k) = H^{0} (\Pp ^* (M^{2n}), (d_k)^*_\om) = H^0 (\Pp ^*(M^{2n}), d_k ^-)= \R$.
\end{corollary}
 
\begin{proof} Note that $\Pp^0 (M^{2n}) = \Om ^0 (M^{2n})$ and $d ^+_k = d_k$, which implies the first  assertion of Proposition
\ref{lich0} immediately. The second assertion of Proposition \ref{lich0} follows from Proposition \ref{primi} and Proposition \ref{dpm1}, taking into account the equalities $H^{0} (\Pp ^*(M^{2n}), d) = H^0(M^{2n}, \R ) = \R$.
\end{proof}

\begin{remark} \label{his1} 1.  Let $(M^{2n}, \om, \theta)$ be a l.c.s. manifold. The symplectic star operator $*_\om$ provides  an isomorphism  between the space  $\Pp^*M^{2n}$ of primitive forms and  the space $\Cc ^* M^{2n}$ of
coeffective forms. Thus  $H^{*}(\Pp^* (M^{2n})) $ is isomorphic to  $H^*(\Cc ^* M^{2n}, d_\theta)$.   The latter cohomology groups for symplectic manifolds
have been introduced by Bouche   \cite{Bouche1990}.   A variant  of the effective cohomology groups  for contact manifolds has been introduced  (and computed) by Lychagin  \cite{Lychagin1979}  already in 1979. Later, a modification of this complex for contact manifolds   has been considered by Rumin independently \cite{Rumin1994}.    Chinea, Marrero and Leo  generalized  the construction of effective cohomology groups for Jacobi manifolds \cite{CML1998}.

2. Note that the groups $H^{q} (\Pp^*(M^{2n}),d^+_k)$  have
the following simple interpretation.  We consider  the  differential ideal $  L (\Om ^* (M^{2n})) \subset \Om^*(M^{2n})$.  The quotient  $\Om^* (M^{2n})/ L(\Om^*(M^{2n}))$ is isomorphic to the space $\Pp^*(M^{2n})$,
and  the differential $d_k$ induces the  differential $d_k ^+$ on the quotient complex.  
  
3. The plus-primitive cohomology groups   and the minus-primitive cohomology groups for symplectic manifolds  
have been   introduced by Tseng and Yau  \cite[part II]{TY2009}.  

4. Below we shall show a deeper  relation between these new cohomology groups  and the twisted cohomology groups $H^*(M^{2n}, d_k)$   using the spectral sequence introduced in the next section.

5. Let $(M^{2n+1}, \alpha)$ be a contact manifold.   Then its symplectization $M^{2n+2}:= M^{2n+1} \times \R$ is supplied with a symplectic form $\om(x, t) = \exp ^t (d\alpha + dt \wedge \alpha) = \tilde \alpha$. 
Denote by $i : M^{2n+1}  \to M^{2n+2}$ the embedding $x \mapsto (x, 0)$.  We observe that the restriction of the 
filtration on $(M^{2n +2}, d \tilde \alpha)$  to $i (M^{2n+1})$ coincides with the    filtration   introduced
by Lychagin \cite{Lychagin1979}. 

6. Note that any symplectic manifold $(M^{2n}, \om)$ is a Poisson manifold equipped with the Poisson bivector $G_\om$.  Associated with  a Poisson bivector $\Lambda$ on a Poisson manifold
$M$ there are two differential complexes. The first one is  the complex of multivector fields on $M$ equipped with the differential $\Lambda$ acting via the Schouten bracket. It has been introduced by  Lichnerowicz and the associated cohomology group is called
the Lichnerowicz-Poisson cohomology of $M$ \cite{Lichnerowicz1977}, \cite{FIL1996}.  The second  differential complex is the complex of differential forms on $M$ equipped with the differential $\delta = [i(\Lambda), d]$ where
$i(\Lambda)$ is the contraction with $\Lambda$.  This complex has been introduced by Kozsul in \cite{Koszul1985} and it is called  the canonical Poisson  homology of $M$ \cite{Brylinski1988}.
 If $M^{2n}$ is symplectic then  $G_\om\in End (T^*M^{2n}, TM^{2n})$  induces an isomorphism between  the de Rham cohomology and  Lichnerowicz-Poisson cohomology \cite[Theorem 6.1]{FIL1996}, and the symplectic star operator   provides an isomorphism between the de Rham cohomology and the canonical Poisson homology \cite{Brylinski1988}.  In \cite{FIL1996} the authors consider the coeffective  Lichnerowicz-Poisson cohomology groups on a Poisson manifold, which are isomorphic to the coeffective symplectic groups introduced by Bouche \cite{Bouche1990} if the Poisson structure is symplectic.
 
 \end{remark}

\section{Spectral sequences  on a l.c.s. manifold}

In this section,   for any integer $k$, we  construct   a   spectral sequence associated with the filtered complex $(F^*K^*, d_k)$   on a l.c.s. manifold  $(M^{2m}, \om, \theta)$.   We compare the   $E_1$-term  of this spectral sequence  with the primitive   cohomology group   $H^{*,+}(\Pp^*(M^{2n}), d_k)$ introduced in the previous section (Lemma \ref{fk1}).
We show the existence of a long exact sequence  connecting  the $E_1$-term of this spectral sequence  with  the  Lichnerowicz-Novikov cohomology
groups $H^* (\Om^*(M^{2n}), d_{k+p})$ for appropriate $p$, which will be denoted by $H^* _{k+p}(M^{2n})$ (Theorem \ref{ex1}, Proposition \ref{prol1}).   We prove that the $E_1$-term of our spectral sequence
is a conformal invariant of $(M^{2n},\om, \theta)$, moreover, the $E_1$-terms associated with $(M^{2n}, \om, \theta)$
and $(M^{2n}, \om', \theta)$ are  isomorphic, if $\om'= \om + d_{\theta} \tau$ (Theorem \ref{conf}).

 In  Proposition \ref{clo2} of the previous section we proved that  $(F^*K^*, d_k)$ is a filtered  complex.
Let us study  the  spectral sequence  $\{ E_{k,r}^{p,q}, d_{k,r}: E^{p,q}_{k,r} \to 
E^{p+r,q-r +1}_{k,r}\}$    of this filtered complex, first introduced by Di Pietro and Vinogradov for  symplectic manifolds in \cite{PV2006}. We refer the reader to  \cite{McCleary2001}, \cite{GH1978}  for an introduction
into the theory of spectral sequences associated with a filtration.
The initial term  $E^{p,q}_{k,0}$ of  this  spectral sequence   is defined as follows, 
\begin{equation}
E^{p,q}_{k,0} = F^p K^{p+q}/F^{p+1} K^{p+q}.\label{e00}
\end{equation}
 
Using the induced Lefschetz decomposition  (\ref{fil1}), (\ref{fil1a}), 
taking into account the injectivity of the map $L ^p : \Om ^{q-p}(M^{2n}) \to \Om^{q+p} (M^{2n})$, we get  for all $k \in \Z$
\begin{equation}
E^{p,q}_{k,0} \cong \Ll^{p, q-p} \cong \Ll ^{0, q-p} \text{ if }  n \ge  q   \ge p \ge 0,
\label{e01}
\end{equation}
\begin{equation}
E^{p,q}_{k,0} = 0 \text{ otherwise }.
\label{e02}
\end{equation}

\medskip

Since $E^{p,q}_{k,l}$ is a quotient of $E^{p,q}_{k, l-1}$, in view of (\ref{e02}), any term $E^{p,q}_{k,l}$  written below, if without explicit condition on $p$ and $ q$,   is   always under the assumption $ 0\le p \le q \le n$.

Now let us go to the next term $E^{p,q}_{k,1}$ of our spectral sequence.
Recall that  $d_{k, 0} :  E^{p,q}_{k,0} \to E^{p, q+1}_{k,0} $ is obtained  from the  differential $d_{k}$  by passing to the quotient:

\begin{equation}
\xymatrix{E^{p,q}_{k,0}   \ar[r]^{d_{k,0}}\ar@{=}[d] &   E^{p, q+1}_{k,0}\ar@{=}[d] \\
F^p K^{p+q}/F^{p+1}K^{p+q}\ar[r] & F^p K^{p+q+1}/F^{p+1}K^{p+q+1} & 
}
\label{e0d}
\end{equation}

Let us write  $d_{k,0}$   explicitly using  (\ref{e01}) and (\ref{e02}). 
Since $d_k L^p = L^p d_{k-p}$, using (\ref{e00}), (\ref{e01}), (\ref{e02})  and (\ref{hoya}) we have for any $\alpha \in  E^{p,q}_{k,0}$ with  $n \ge q \ge p \ge 0$
\begin{equation}
d_{k,0}  (\alpha)  =[L^p (d_{k-p} ^+ + L d_{k-p} ^-)( \tilde \alpha)] =  [L^p d_{k-p} ^ + ( \tilde \alpha)]\in E^{p, q+1}_{k,0}, \label{vy2}
\end{equation}
where $\tilde \alpha \in \Ll^{0, q-p}$  is  a representative of $\alpha \in  E^{p,q}_{k,0}$ by (\ref{e01}).
Since $L^p : \Ll ^{0, q-p} \to \Ll^{p, q-p}$ is  an isomorphism, if $0\le p \le q \le n$ by (\ref{e01}),  we rewrite  (\ref{vy2}) as follows
\begin{eqnarray}
d_{k,0}: E^{p,q}_{k,l}= \Ll ^{0, q-p} \to E^{p, q+1}_{k, 0} = \Ll ^{0, q+1 -p}_{ k, 0}, \, \tilde \alpha \mapsto d^+_{k-p}\tilde \alpha, \label{e12} \\
\text{  if }   0\le p \le q \le n. \nonumber 
 \end{eqnarray}

\begin{lemma}\label{fk1} The  term $E^{*,*}_{k,1}$ of the  spectral sequence $\{ E^{p,q}_{k,r}, d_{k,r}\}$ is  determined by the following relations
\begin{equation}
E^{p,q}_{k,1} = H^{q-p}  (\Pp^*(M^{2n}), d^+_{k-p}) \text { if } 0 \le p \le q \le n-1,
\label{e13}
\end{equation}
\begin{equation}
E^{p,n}_{k,1} = { \Pp^{n-p}(M^{2n}) \over  d_{k-p}^+ (\Pp^{n-p-1}(M^{2n}))}, \text { if } 0 \le p \le n,\label{e13a}
\end{equation}
\begin{equation}
E^{p,q}_{k,1} = 0 \text { otherwise }. \label{e13b}
\end{equation}
\end{lemma}

\begin{proof}  Clearly (\ref{e13}) is a consequence of (\ref{e12}).
Next using (\ref{vy2}) and the  identity $L ^p( \Pp^{n -p +1 } (M^{2n})) = 0$, we obtain  $d_{k, 0} (E^{p, n} _{k,0} ) = 0$.   Hence
\begin{equation}
E^{p,n}_{k,1} = { L^p(\Pp^{n-p}(M^{2n})) \over L^p ( d_{k-p}^+ (\Pp^{n-p-1}(M^{2n})))}. \label{e13aa}
\end{equation}
Since $L^p : \Om ^{n-p} (M^{2n}) \to \Om ^{n+p} (M^{2n})$ is injective, (\ref{e13aa}) implies (\ref{e13a}).

The last  relation (\ref{e13b}) in Lemma \ref{fk1} follows from (\ref{e02}). This completes the proof of Lemma \ref{fk1}.
\end{proof}

Next we  define the following diagram  of  chain complexes
$$
\xymatrix{ \Om ^{q-p-1}(M^{2n})\ar[r] ^ {L} \ar[d] ^{d_{l-1}} & \Om ^{q-p+1}(M^{2n}) \ar[r]^{\Pi L^p}\ar[d]^{d_{l} }   & E^{p, q+1}_{l+p,0}\ar[d]_{d_{l+p,0}}\ar[r]&\\
\Om ^{q-p}(M^{2n}) \ar[r]  ^L\ar[d]^{d_{l-1}}  &  \Om ^{q-p+2}(M^{2n}) \ar[r]^ {\Pi L^p}\ar[d]^{d_{l}}   &E^{p, q+2}_{l+p,0} \ar[d]_{d_{l+p,0}}\ar[r] & \\
\Om ^{q-p+1}(M^{2n})\ar[r] ^L   &  \Om ^{q-p+3}(M^{2n}) \ar[r]^{\Pi L^p} & E^{p, q+3}_{l+p,0} \ar[r] &
}
$$

Here  the map $\Pi L^p$ associates  with each element  $\beta \in \Om^{q-p+1}(M^{2n})$  the element $[L ^p\beta]\in
E^{p, q+1}_{l+p,0}$. 
Recall that   $d_{l} \circ L = L \circ d_{l-1}$.  Thus the upper part of the  above  diagram is commutative.
The lower part of the diagram is also commutative, since
$$  d_{l+p} L^ p  = L^p d_{l} .$$

Summarizing we have the following  sequence  of  chain complexes
\begin{equation}
0\to (\Om^{q-(p+1)}(M^{2n}), d_{l-1})  \stackrel{L}{\to} (\Om^{q+1 -p}(M^{2n}), d_{l}) \stackrel{\Pi L^p}{\to} (E_{l+p,0}^{p, q+1}, d_{l+p,0}) \to 0.\label{exac3}
\end{equation}

Set $\Om ^{-1} (M^{2n}): = 0$.

\begin{lemma}\label{lexac3}  The  sequence (\ref{exac3}) of chain complexes  is exact  for $0\le p \le q \le n-1$.
\end{lemma}

\begin{proof}  
 For $0\le p \le q  \le n$ the operator $L: \Om ^{q-p-1} (M^{2n}) \to \Om ^{q-p+1}(M^{2n})$ is injective by Lemma \ref{inj1}, so the  sequence
(\ref{exac3}) is exact  at $\Om ^{q -(p+1)}(M^{2n})$. The exactness at $\Om^{q+1-p}(M^{2n})$ follows easily from Corollary \ref{inj1}, taking into account (\ref{e01}).  The exactness at $E^{p,q+1}_{l+p,0}$ follows directly from the definition (\ref{e00}) of $E^{p,q}_{l+p,0}$. 
\end{proof}

As a consequence of Lemma \ref{lexac3}, using the general theory of homological  algebra, see e.g. \cite[Chapter 1, \S 6]{GM1988}, we get immediately 

\begin{theorem}\label{ex1} The following long  sequence is exact  for $0\le p \le q \le n-1$
\begin{equation}
\cdots \to H^{q-p}_l(M^{2n}) \stackrel{ \bar L^p}{\to} E^{p,q}_{l+p,1}\stackrel{\delta_{p,q}}{\to} H^{q-(p+1)}_{l-1}(M^{2n})\stackrel{\bar L} {\to} H^{q+1-p}_{l}(M^{2n}) \stackrel{\bar L^p} {\to } E^{p, q+1}_{ l+p,1} \to \cdots \label{ex11}
\end{equation}
\end{theorem}

\begin{remark}\label{rpv} 1. Theorem \ref{ex1} is  a generalization of  \cite[Theorem 1]{PV2006} stated for symplectic manifolds.

2. Let us write the connecting homomorphism $\delta_{p,q}$  explicitly.    If $\alpha \in H^{q-p} (\Pp^*(M^{2n}), d_l^+) =
E^{p,q}_{l+p,1}$, see (\ref{e13}), then $d_l\tilde \alpha = L \eta $   for any representative
$\tilde \alpha \in \Pp^{q-p}(M^{2n})$ of $\alpha$. By the definition of connecting homomorphism, see also \cite[\S 3]{PV2006}, $\delta_{p,q} (\alpha) = [\eta ] \in H^{q-(p+1)}_{l-1} (M^{2n})$. Since $d^+_l\tilde \alpha = 0$, $d_l \tilde \alpha = Ld_l ^- \tilde \alpha$, so $\eta = d^-_l \tilde \alpha$. Thus we get
\begin{equation}
\delta_{p,q}  \alpha = [d_{l}^- \tilde \alpha] \in H^{q-(p+1)}_{l-1} (M^{2n}).\label{delta1}
\end{equation}

3. Let us write the operator $\bar L^p$ explicitly.  Assume that $[\beta] \in H^{q-p} _l (M^{2n})$, $\beta \in \Om^{q-p} (M^{2n})$. Set 
$$\beta _{pr}: = \Pi_{pr} \beta$$
 - the primitive component  of $\beta$ in the first Lefschetz decomposition.  Since $d_{l}\beta = 0$, we have $d_{l} ^ + \beta _{pr} = 0$.
Thus the image  $\bar L^p [\beta] = [L^p \beta] \in E^{p,q}_{l+p, 1}=H^{q-p} (\Pp^*(M^{2n}), d_l^+)$ has a  representative $\beta _{pr}\in \Pp^{q-p}(M^{2n})$ with 
$d_l^+\beta _{pr} = 0$. Summarizing we have
\begin{equation}
\bar L^p [\beta] = [\beta_{pr}] \in H^{q-p} (\Pp^*(M^{2n}), d^+_l) = E^{p,q}_{l+p, 1}. \label{lpex}
\end{equation}
\end{remark}

4. Substituting for $0\le p= q \le n-1$ in the  long exact sequence (\ref{ex11}), we get  the left end of (\ref{ex11})
\begin{equation}
0 \to  H^0_l (M^{2n}) \to E^{p,p}_{l+p, 1}  \to 0  \to H^1 _l (M^{2n}) \to \cdots .  \label{dege1} 
\end{equation}
From  (\ref{e13})  and (\ref{dege1}) we get for $0\le p \le n-1$
\begin{equation}
E^{p, p}_{l+p, 1}= H^{0,+} (\Pp^*(M^{2n})) = H ^0_l (M^{2n}) .\label{0pl}
\end{equation} 
\medskip
Let us  prolong the exact sequence (\ref{ex11}) for $q=n$, using the ideas in \cite[\S III]{PV2006}. 
 For $0 \le k \le n$ we set
 $$C^k_l : = {\ker d_l ^- \cap \Om ^k (M^{2n})\over  d_l ( \Om ^{k-1} (M^{2n}))}.$$

\begin{proposition}\label{prol1} The long exact sequence (\ref{ex11})  can be extended as follows
\begin{equation}
E^{p, n-1}_{l+p, 1}\stackrel{\delta_{p,n-1}}{ \to} H^{n - (p+2)}_{l-1}(M^{2n}) \stackrel{[L]}{\to} C^{n-p}_{l}\stackrel{[ L^p] }{\to} E^{p, n}_{l+p,1} \stackrel{\delta_{p,n}}{\to}  C^{n - (p+1)}_{l-1} \stackrel{[ L^{p+1}] }{\to} H^{n +p +1}_{l+p} (M^{2n}) \to 0,\label{prol1a}
\end{equation}
where $0\le p \le n-1$ and the operators $[L], [L^p]$ and $[L^{p+1}]$ will be defined in the proof  below.
\end{proposition}

\begin{proof} First we define $[L]$ and prove the exactness at $H^{n-(p+2)}_{l-1} ( M^{2n})$. Denote by $\Pi: H^{n-p}_{l} (M^{2n}) \to  C^{n-p}_{l}$ the natural embedding of the quotient spaces
$$\frac{\ker d_l \cap \Om ^{n -p} (M^{2n})}{d_l (\Om ^{n-p-1} (M^{2n}))} \to \frac{\ker d_l^- \cap \Om ^{n -p} (M^{2n})}{d_l (\Om ^{n-p-1} (M^{2n}))}.$$
 Set $[L]: = \Pi \circ \bar L$, where $\bar L:H^{n - (p+2)}_{l-1}(M^{2n}) \to  H^{n-p}_{l}(M^{2n})$ is induced by $L$. By Theorem \ref{ex1} we
have  $Im \, \delta _{p, n-1} = \ker \bar L$.
Since $\Pi$ is an embedding, the last  equality implies $\ker [ L ]= \ker \bar L$.  This proves the  required exactness at $H^{n-(p+2)}_{l-1} ( M^{2n})$.

Now we  define $[L^p]$ and show the exactness at $C^{n-p}_{l}$.  Assume that $\alpha = \alpha _{pr} + L \tilde \alpha \in \Om ^{n-p} (M^{2n})$ is a representative of $[\alpha]\in  C^{n-p}_{l}$,
i.e. $ d_{l}^- (\alpha ) = 0$, or equivalently $d_l \alpha = d_l ^+ \alpha$. 
We  set
\begin{equation}
[L^p] (\alpha): = [\alpha_{pr}]\in  \frac{\Pp^{n-p}(M^{2n})}{ d _{l}^+ (\Pp^{n-p-1}( M^{2n}))} = E^{p,n} _{l+p, 1} .\label{new1a}
\end{equation}
Clearly  the map $[L^p]$ is well-defined, since $\Pi_{pr}d_l (\gamma) = d_l ^+ \Pi_{pr}(\gamma)$.
Now assume  that $[\alpha] \in \ker [L^p]$, so  by (\ref{new1})
\begin{equation}
\Pi_{pr}\alpha = d^+_l  \gamma \text  { for  some }\gamma \in \Pp ^{n-p-1} (M^{2n}).\label{new11}
\end{equation} 
Using  the property $d_l \alpha  = d_l ^+ \alpha$ we obtain   from (\ref{new11}) $d_l \alpha = 0$.
Now we write 
\begin{equation}
\alpha = \alpha _{pr} + L\tilde \alpha = d_l ^+ \gamma + L\tilde \alpha =  d_l \gamma + L \beta,\label{new2}
\end{equation}
where $\beta = d_l^- \gamma + \tilde \alpha$.
Since $d_l \alpha = 0$, using (\ref{new2}) we get  $d_l L\beta = L d_{l-1} \beta = 0$.  Applying Lemma \ref{inj1} to  $d_{l-1}\beta \in \Om^{n-p-1}(M^{2n})$, we obtain $d_{l-1} \beta = 0$. This implies $[\alpha] = [L] ([\beta]) \in   \im\, [L]$, and the required exactness. 

Next we define the operator $\delta_{p,n}:E^{p, n}_{l+p,1} \to  C^{n - (p+1)}_{l-1} $ as follows
\begin{equation}
\delta_{p,n}(\alpha ) :=  [d_l ^- \tilde \alpha] \in C^{n - (p+1)}_{l-1},\label{deltac}
\end{equation}
where $\tilde \alpha \in \Pp^{n-p}( M^{2n})$ is a representative of  $\alpha \in E^{p,n} _{l+p,1}= \Pp^{n-p}(M^{2n})/ d _{l}^+ (\Pp^{n-p-1} (M^{2n}))$.
Clearly $\delta_{p,n} (\alpha) \in C^{n-p-1}_{l-1}$, since by (\ref{com1}) $d_{l-1} (d_l^-\tilde  \alpha ) = d^+_{l-1} d^-_{l } \tilde \alpha$.
 The map $\delta_{p,n}$  is well-defined, since   for any $\gamma \in \Pp^{n-p-1}(M^{2n})$  using (\ref{com1}) and (\ref{com1a}) we get
\begin{equation}
[d_l ^- d_l^+ \gamma] = -[d_{l-1}^+ d_l ^- \gamma]  = -[d_{l-1} (d_l^- \gamma)] = 0 \in C^{n-p-1}_{l-1}.
\end{equation}
Now assume that $\alpha \in \ker \delta _{p,n}$ and $\tilde \alpha \in \Pp^{n-p}(M^{2n})$ is its representative. By (\ref{deltac})  $d_{l} ^- \tilde \alpha = d_{l-1} \beta$ for some
$\beta \in \Om ^{n-p-2} (M^{2n})$. It follows
\begin{equation}
d_{l} ^+ \tilde \alpha = d_l  \tilde \alpha -  L d_{l-1} \beta =  d_l ( \tilde \alpha - L\beta).\label{deltac1}
\end{equation}
Since $\tilde \alpha$ is primitive,  and $d^+_l (\Ll ^{s, r}) \subset \Ll ^{s, r+1}$, we get
\begin{equation}
 d_{l}^+ \tilde \alpha =  d_l ^+ (\tilde \alpha -  L \beta).\label{deltac2}
\end{equation}

Clearly, (\ref{deltac1}) and (\ref{deltac2}) imply  $[\tilde \alpha -  L  \beta] \in C^{n-p}_l$,  and   by
(\ref{new1})  $\alpha = [L^p] ([\tilde \alpha -  L  \beta])$. This yields the exactness at $E^{p, n}_{l+p,1}$. 

Let us define $[L^{p+1}]$ and  show the exactness at $C^{n-(p+1)}_{l-1}$.  For $\alpha \in C^{n-p-1}_{l-1}$  we set
\begin{equation}
[L^{p+1}](\alpha): =  [L ^{p+1} \tilde \alpha] \in H^{n+p+1}_{l+p}(M^{2n}),
\end{equation}
 where $\tilde \alpha \in \Om ^{n-p-1}(M^{2n})$ is a representative of $\alpha$.  Note that $d_{l+p} L^{p+1} \tilde \alpha = L^{p+1}d_{l-1}\tilde \alpha = 0$, so $[L^{p+1}] (\alpha) \in H^{n+p+1}_{l+p}(M^{2n})$.  The same formula  shows that our map $[L^{p+1}]$ does not depend on the choice
 of a representative $\tilde \alpha$ of $\alpha$.
Now assume  that $\alpha \in \ker [L^{p+1}]$. Then $L^{p+1}\tilde \alpha = d_{l+p}  \beta$
for some $\beta \in \Om ^{n+p} (M^{2n})$. Using the Lefschetz decomposition for $\beta$ we write
$$\beta = L^p(\beta _{0} + \sum_ {k =1} ^{[{n+p\over 2}]}L^k \beta_{k}), \, \beta_i \in \Pp^*(M^{2n}).$$
It follows that 
\begin{equation}
L^{p+1}\tilde\alpha = d_{l+p}\beta  = L^p (d^+_{l} \beta_{0} + L (d_{l} ^- \beta_{0} + d_{l-1}( \beta_{1} + L\beta _{2} + \cdots))).\label{deltac3}
\end{equation}
By Corollary \ref{inj1}.1
$L^{p+1}: \Om ^{n-p-1}(M^{2n}) \to \Om ^{n+p+1}(M^{2n})$ is an isomorphism, hence we get from (\ref{deltac3})
\begin{equation}
\tilde\alpha = d_{l} ^- \beta_{p} + d_{l-1} (\beta_{p+1} + L\beta_{p+1} + \cdots ). \label{deltac4}
\end{equation}
Combining  (\ref{deltac4}) with  (\ref{deltac}) we get $\alpha \in \im\, \delta_{p,n}$. This proves the required exactness.

Finally we show that $[L^{p+1}]$ is surjective. Assume that $ \beta \in  \Om ^{n+p+1}( M^{2n})$ is a representative of
$[\beta] \in H^{n+p+1}_{l+p+1} (M^{2n})$. Using the Lefschetz decomposition we  write $\beta = L^{p+1} (\tilde \beta)$, $\tilde \beta \in \Om ^{n-p-1}(M^{2n})$. 
Note that $L^{p+2} d_l ^-\tilde \beta = L d_{l+p} ^- \beta = 0$, since $d_{l+p} \beta = 0$.  Since $L^{p+2} : \Om ^{n-p-2}(M^{2n}) \to \Om ^{n+p +2} (M^{2n})$ is an isomorphism, we get $d_l ^- \tilde \beta = 0$. Hence  $[\tilde \beta ] \in C^{n-p-1}_l$ and $[\beta ] =[L^{p+1}]( [\tilde \beta])$.
 This completes the proof of  Proposition \ref{prol1}. 
\end{proof}

\begin{theorem}\label{conf} (cf. \cite[Osservazione 18]{Pietro2006}) The   spectral sequences $E^{p,q}_{k,r}$ on $(M^{2n}, \om, \theta)$ and
on $(M^{2n}, \om ', \theta')$  are isomorphic, if $\om$ and $\om'$ are conformal equivalent. Furthermore, the  $E^{*,*}_{k,1}$-terms  of the spectral sequences on $(M, \om, \theta)$   and $(M, \om', \theta')$ are  isomorphic, if  $\om ' = \om + d_{\theta} \rho$ for some $\rho \in \Om ^1 (M^{2n})$.
\end{theorem}

\begin{proof}  If $\om' = \pm e^f \om$, then  $d_\theta\om ' = df \wedge \om'$. Hence
\begin{equation}
(d_\theta - df\wedge)\om' = 0.\label{confn1}
\end{equation}
Since $L : \Om ^1 (M ^{2n})  \to \Om ^3 (M^{2n})$ is injective,  (\ref{confn1}) implies that 
\begin{equation}
d_{\theta '} = d_{\theta} - df\wedge.\label{confn2}
\end{equation}
It follows
\begin{equation}
d_{k\theta'} = d_{k\theta} - k\cdot df\wedge.\label{confn2}
\end{equation}
 Hence the map $I_f: \alpha \mapsto e^{kf} \alpha$ is an isomorphism between complexes   $(F^* K^*, d_{k\theta})$ and
$(F^*K^*, d_{k\theta'})$, since
\begin{equation}
d_{k\theta'} ( e ^{kf }\alpha) =  d_{k\theta} ( e^{kf} \alpha) + ( - k\cdot df) \wedge e ^{k\cdot f} \alpha =  e^{k\cdot f }(d_{k\theta} \alpha).\label{conf1}
\end{equation}
It follows that the resulting  terms $E^{p,q}_{k,0}$ are also  conformal equivalent.  Moreover, the  map $I_f$ induces an isomorphism $I_f ^0$ between complexes
$$I_f ^0: (E^{p, q}_{k, 0}, d_{k, 0}) \to (E^{p, q}_{k, 0}, d_{k, 0} '), \, [\alpha] \mapsto [e^{kf}\alpha].$$
Inductively, this proves the first assertion of Theorem \ref{conf}. 

 Now we assume that $\om' = \om +d_\theta\rho$.  Then  $d_\theta \om' = 0$.  Using the injectivity of $L : \Om ^1 (M ^{2n})  \to \Om ^3 (M^{2n})$  we conclude that the Lee form of $\om'$  is equal to the Lee form  $\theta$ of $\om$.
Denote by $L'$  the wedge product with $\om'$, and by $[L']$  the induced  operator on $H^*_l (M^{2n})$. Using $\om' - \om = d_\theta \rho$ and applying  (\ref{comd2}), which implies that the wedge product with $d_\theta \rho$ maps $H^k_l (M^{2n})$ to zero, we conclude that the operators $L$ and $L'$   induce the same map $H^k_l (M^{2n}) \to H^{k+2}_{l+1} (M^{2n})$. 

To prove the second assertion of Theorem \ref{conf} we  use the  following version of Five Lemma, whose proof  is obvious and hence omitted.

\begin{lemma}\label{five}
Assume that the following diagram of vector spaces $A_i, B_i$ over a field $\F$ is commutative. If the rows are  exact and $A_1 \to B_1$, $A_2 \to B_2$, $A_4 \to B_4$, $A_5 \to B_5$ are isomorphisms, then there is an  isomorphism from $A_3$ to $B_3$,  which  also commutes with the  other arrows. 
$$
\xymatrix{A_1 \ar[r] \ar[d] & A_2 \ar[r]\ar[d] & A_3 \ar[r] & A_4 \ar[r]\ar[d] & A_5\ar[d]\\
 B_1 \ar[r]  & B_2 \ar[r] & B_3 \ar[r] & A_4 \ar[r] & A_B }
 $$
\end{lemma}
The second assertion of Theorem \ref{conf} for $E^{p,q}_{k,1}$ follows immediately from the  long exact sequence
(\ref{ex11}) and the formula (\ref{conf1}), if $0\le p \le q \le n-1$, taking into account Lemma \ref{fk1}. 


To examine the term $E^{p, n}_{k,1}$, $0\le p \le n-1$, we need the following

\begin{lemma}\label{five2} 
Assume that $\omega'=\omega+d_\theta\rho$. For $0\le p \le n$ there are linear maps $B^{n-p}_l:C^{n-p}_l(\omega)\rightarrow
C^{n-p}_l(\omega')$  such that the following two diagrams are commutative. (The symbol $I$ denotes the identity
mapping. The other mapping are defined in the proof.)

\xymatrix{
H^{n-(p+2)}_{l-1}(M^{2n})\ar[d]^I\ar[r]^{[L]} & C^{n-p}_l(\omega)\ar[d]^{B^{n-p}_l} &
C^{n-(p+1)}_{l-1}(\omega)\ar[d]^{B^{n-p-1}_{l-1}}\ar[r]^{[L^{p+1}]} & H^{n+p+1}_{l+p}(M^{2n})\ar[d]^I\\
H^{n-(p+2)}_{l-1}(M^{2n})\ar[r]^{[L']} & C^{n-p}_l(\omega') &
C^{n-(p+1)}_{l-1}(\omega')\ar[r]^{[L^{\prime p+1}]} & H^{n+p+1}_{l+p}(M^{2n})
}
\end{lemma}
\begin{proof}
Let us define first a linear mapping 
$$
\tilde{B}^{n-p}_l:\ker d^-_l(\omega)\cap\Omega^{n-p}(M^{2n})\rightarrow\Omega^{n-p}(M^{2n}).
$$
Let $\eta\in\ker d^-_l\cap\Omega^{n-p}(M^{2n})$. This means that $d_l\eta=d_l^+\eta$ or equivalently 
that $d_l\eta$ is primitive. The last assertion is again equivalent to the equality $\omega^p\wedge d_l\eta=0$. Since $(L')^p : \Om ^{n-p}(M^{2n}) \to \Om ^{n+p} (M^{2n})$ is an isomorphism,  there is a unique $\eta'\in\Omega^{n-p}(M^{2n})$ such that
$$
\sum_{i=1}^p\binom{p}{i}\rho\wedge(d_\theta\rho)^{i-1}\wedge\omega^{p-i}\wedge d_l\eta=\omega^{\prime p}\wedge\eta'.
$$
Now we define $\tilde B^{n-p}_l$ by
$$\tilde{B}^{n-p}_l\eta:=\eta-\eta'.$$
 We shall now prove that the element $d_l(\eta-\eta')$ is $\omega'$-primitive.
\begin{gather}
\omega^{\prime p}\wedge d_l(\eta-\eta')=-\omega^{\prime p}\wedge d_l\eta'+(\omega+d_1\rho)^p\wedge d_l\eta=\notag\\
=-d_{p+l}(\omega^{\prime p}\wedge\eta')+\omega^p\wedge d_l\eta+
\sum_{i=1}^p\binom{p}{i}(d_1\rho)^i\wedge\omega^{p-i}\wedge d_l\eta=\notag\\
=-d_{p+l}(\omega^{\prime
p}\wedge\eta')+d_{p+l}\bigg[\sum_{i=1}^p\binom{p}{i}\rho\wedge(d_1\rho)^{i-1}\wedge\omega^{p-i}
\wedge d_l\eta\bigg]=\notag\\
=d_{p+l}\bigg[-\omega^{\prime p}\wedge\eta'+\sum_{i=1}^p\binom{p}{i}\rho\wedge(d_{1}\rho)^{i-1}\wedge\omega^{p-i}
\wedge d_l\eta\bigg]=0.\notag
\end{gather}
In other words we have proved that $\tilde{B}^{n-p}_l$ maps $\ker d^-_l\cap\Omega^{n-p}(M^{2n})$ into $\ker (d')^-_l\cap\Omega^{n-p}(M^{2n})$, where $(d')^-_l$ is defined via  the Lefschetz decomposition corresponding to $\om'$.

Let us take now an element $\eta\in d_l\Omega^{n-p-1}(M^{2n})$, i. e. $\eta=d_l\gamma$, where $\gamma\in \Omega^{n-p-1}(M^{2n})$. Then
$d_l\eta=0$ and we have $\omega^{\prime p}\wedge\eta'=0$, which implies $\eta'=0$. We thus get $\tilde{B}^{n-p}_l\eta=\eta$.
Consequently we have proved that $\tilde{B}^{n-p}_l$ maps $d_l(\Omega^{n-p-1}(M^{2n}))$ into itself. Now it is obvious that $\tilde{B}^{n-p}_l$
induces a linear mapping 
$$B^{n-p}_l:C^{n-p}_l(\omega)\rightarrow C^{n-p}_l(\omega').$$

Next we shall investigate the first diagram. First we define the mapping $[L]$. If $[\beta]\in H^{n-(p+2)}_{l-1}
(M^{2n})$, then we have an element $\beta\in\Omega^{n-(p+2)}(M^{2n})$ such that $d_{l-1}\beta=0$. Let us set 
$$[L][\beta]:=[\omega\wedge\beta].$$
 It is easy to see that this definition depends only on the cohomology class $[\beta]\in H^{n-p-2}_{l-1} (M^{2n})$. Namely, if $\tilde{\beta}=\beta+d_{l-1}\gamma$, then
$$
\omega\wedge(\beta+d_{l-1}\gamma)=\omega\wedge\beta+d_l(\omega\wedge\gamma),
$$
which shows that $[L][\tilde{\beta}]=[L][\beta]$. Similarly we define $[L']$. Let us take $[\beta]\in H^{n-(p+2)}_{l-1}
(M^{2n})$. Then $d_{l-1}\beta=0$, and we have
\begin{gather}
\sum_{i=1}^p\binom{p}{i}\rho\wedge(d_1 \rho)^{i-1}\wedge\omega^{p-i}\wedge d_l(\omega\wedge\beta)=\notag\\
=\sum_{i=1}^p\binom{p}{i}\rho\wedge(d_1 \rho)^{i-1}\wedge\omega^{p-i}\wedge(d_1\omega\wedge\beta+\omega\wedge
d_{l-1}\beta)=0.\notag
\end{gather}
We have $0=\omega^{\prime p}\wedge\eta'$, which implies $\eta'=0$. Thus we get $B^{n-p}_l[L][\beta]=
B^{n-p}_l[\omega\wedge\beta]=[\omega\wedge\beta-0]=[\omega\wedge\beta]$.
On the other hand we  compute
\begin{gather}
[L'][\beta]=[\omega'\wedge\beta]=[(\omega+d_1\rho)\wedge\beta)]=[\omega\wedge\beta]+[d_1\rho\wedge\beta+\rho\wedge
d_{l-1}\beta]=\notag\\
=[\omega\wedge\beta]+[d_l(\rho\wedge\beta)]=[\omega\wedge\beta].\notag
\end{gather}
We have thus shown that that the first diagram is commutative.

We continue now with the second diagram. Again we define first the mapping $[L^{p+1}]$. For $[\eta]\in
C^{n-(p+1)}_{l-1}(\omega)$ there is a representative $\eta\in\Omega^{n-(p+1)}_{l-1}(M^{2n})$ such that $d^-_{l-1}\eta=0$. 
This means that $d_{l-1}=d^+_{l-1}\eta$. We compute
\begin{gather}
d_{l+p}(\omega^{p+1}\wedge\eta)=d_{(p+1)+(l-1)}(\omega^{p+1}\wedge\eta)=d_{p+1}(\omega^{p+1})\wedge\eta+\omega^{p+1}
\wedge d_{l-1}\eta=\notag\\
=\omega^{p+1}\wedge d^+_{l-1}\eta=0.\notag
\end{gather}
The last term vanishes because the form $d^+_{l-1}\eta$ is primitive. Finally let us suppose that
$\eta=d_{l-1}\gamma$. Then we have
$$
\omega^{p+1}\wedge d_{l-1}\gamma=d_{p+1}(\omega^{p+1})\wedge\gamma+\omega^{p+1}\wedge d_{l-1}\gamma=d_{(p+1)+(l-1)}
(\omega^{p+1}\wedge\gamma).
$$
This shows that we can define $[L^{p+1}]$ by the formula 
$$[L^{p+1}][\eta]:=[\omega^{p+1}\wedge\eta].$$
 Now we are going
to prove the commutativity of the second diagram. For $[\eta]\in C^{n-(p+1)}_{l-1}(\omega)$ we have $B^{n-p-1}_{l-1}[\eta]=
[\eta-\eta']$, where $\eta'$ is uniquely determined by the equality
$$
\omega^{\prime p+1}\wedge\eta'=\sum_{i=1}^{p+1}\binom{p+1}{i}\rho\wedge(d_1 \rho)^{i-1}\wedge\omega^{p+1-i}\wedge
d_{l-1}\eta.
$$ 
Further we have $[L']B^{n-p-1}_{l-1}[\eta]=[\omega'\wedge(\eta-\eta')]$. Now let us compute
\begin{gather}
\omega^{\prime p+1}\wedge(\eta-\eta')-\omega^{p+1}\wedge\eta=\omega^{\prime p+1}\wedge\eta 
-\omega^{\prime p+1}\wedge\eta'-\omega^{p+1}\wedge\eta=\notag\\
=(\omega+d_1 \rho))^{p+1}\wedge\eta-\sum_{i=1}^{p+1}\binom{p+1}{i}\rho\wedge(d_1 \rho)^{i-1}\wedge\omega^{p+1-i}\wedge
d_{l-1}\eta-\omega^{p+1}\wedge\eta=\notag\\
=\sum_{i=1}^{p+1}\binom{p+1}{i}(d_1 \rho)^i\wedge\omega^{p+1-i}\wedge\eta
-\sum_{i=1}^{p+1}\binom{p+1}{i}\rho\wedge(d_1 \rho)^{i-1}\wedge\omega^{p+1-i}\wedge d_{l-1}\eta=\notag\\
=d_{p+1}\bigg(\sum_{i=1}^{p+1}\binom{p+1}{i}\rho\wedge(d_1 \rho)^{i-1}\wedge\omega^{p+1-i}\bigg)\wedge\eta\notag\\
-\bigg(\sum_{i=1}^{p+1}\binom{p+1}{i}\rho\wedge(d_1 \rho)^{i-1}\wedge\omega^{p+1-i}\bigg)\wedge d_{l-1}\eta=\notag\\
=d_{(p+1)+(l-1)}\bigg(\sum_{i=1}^{p+1}\binom{p+1}{i}\rho\wedge(d_1 \rho)^{i-1}\wedge\omega^{p+1-i}\wedge\eta\bigg).
\end{gather}
We have thus proved that $[L^{\prime p+1}]B^{n-p-1}_{l-1}[\eta]=[L^{p+1}][\eta]$.
\end{proof}

Clearly the second assertion of Proposition \ref{conf} for $E^{p,n}_{k,1}$, $0\le p \le n-1$,  follows from Proposition \ref{prol1},   Lemma \ref{five},   and Lemma \ref{five2}.
Combining with  (\ref{e13a}), which implies that $E^{n,n}_{k,1}  = C^\infty (M^{2n})$, we obtain  the second assertion of Proposition \ref{conf}.
This completes the proof of Theorem \ref{conf}.
\end{proof}

\begin{remark}\label{pois}
In \cite[Example 7.1]{FIL1998} the authors construct an example of a compact 6-dimensional nilmanifold $M^6$ equipped with a family of  symplectic forms $\om_t$, $t \in [0,1],$ with varying cohomology classes $[\om_t] \in H^2 (M^6, \R)$.  They showed that the  coeffective cohomology groups associated to $\om_1$ and $\om_2$ have different Betti number $b_4$. It follows that, using \cite[Lemma 2.7,  Proposition 3.5, part II]{TY2009}, see also Remark \ref{his1}.1 above, 
the $E_1$-terms  of the  associated spectral sequences for $\om_0$ and $\om_1$ are different.
\end{remark}

\section{The stabilization of  the  spectral sequences}

In this section we prove  that the spectral sequences  $\{E^{p,q}_{l,r} \}$ on l.c.s.manifolds
$(M^{2n}, \om, \theta)$  converge  to the Lichnerowicz-Novikov cohomology  $H^*_l(M^{2n})$ at  the second term $E^{*,*}_{l,2}$,
or at the $t$-term  $E^{*,*}_{l,t}$  under some cohomological  conditions posed on $\om$  (Theorems \ref{stab1},  \ref{nil},  \ref{stcom}). As a consequence,  we obtain a relation between  the primitive cohomology groups and  the  de Rham cohomology groups of $(M^{2n}, \om)$, if $(M^{2n}, \om)$ is a compact symplectic manifold. This  gives  an answer to a question posed  by Tseng and Yau in \cite{TY2009}, see Remark \ref{tsey}.

First we   prove   the following simple property of the second  terms $E^{*,*}_{l,2}$, which will be used later in the proof of Theorem \ref{stcom}.

\begin{proposition}\label{symm} (cf. \cite[Proposizione 19]{Pietro2006}) Assume that $1\le p\le q \le n-1$. Then $E^{p, q}_{l, 2} = E^{p-1, q-1}_{l, 2}$.
\end{proposition}

\begin{proof} Let $\alpha \in  E^{p, q}_{l+p, 1} =H^{q-p} (\Pp^* (M^{2n}), d^+_l) $ and $\tilde \alpha \in \Pp^{q-p} (M^{2n})$ its representative  as in (\ref{e13}).  The differential $d_{l+p,1}: E^{p,q}_{l+p, 1} \to E^{p+1, q}_{l+p, 1}$ is defined by
\begin{equation}
d_{l+p,1} (\alpha): = [d_{l+p} L^p \tilde \alpha] = [L^p d_{l} \tilde \alpha] \in E^{p+1, q}_{l+p, 1}.\label{stab1f}
\end{equation}
Using  $d _{l+p} = d_{l+p}^+ + Ld_{l+p}^-$  and taking into account $d^+_{l} \tilde \alpha = 0$,  we observe that
 $[L^p d_l\tilde \alpha] \in E^{p+1, q}_{l+p, 1}$ has a representative $d_{l}^- (\tilde \alpha) \in \Pp^{q-p-1}  (M^{2n})$    in  $ H^{q-p-1} (\Pp^*( M^{2n}), d^+_{l-1}) =E^{p+1, q}_{l+p, 1} $, if $0\le p \le q \le n$. Equivalently, using  (\ref{e13}),  we rewrite  $d_{l+p,1}$   for $0\le p \le q \le n$ as follows
\begin{equation}
d_{l+p,1} : H^{q-p} (\Pp ^*(M^{2n}), d^+_l) \to  H^{q-p-1} (\Pp^*(M^{2n}), d^+_{l-1}), \, [\tilde \alpha] \mapsto [d_l^- \tilde \alpha]. \label{stab1a}
\end{equation}
Clearly Proposition \ref{symm}   follows from  (\ref{e13}) and the formula (\ref{stab1f}), (\ref{stab1a}).
\end{proof}

Now assume that  $\om = d_k\tau$ for some $k \in \Z$  and $\tau \in \Om^1(M^{2n})$. Since $d_1 \om =  d_k \om = 0$, it follows that $(k-1)\theta \wedge \om =0$. Since $L$ is injective, we get $k = 1$.   The following  theorem is a generalization of \cite[Theorem 2]{PV2006} for the symplectic case $\theta = 0$.

\begin{theorem}\label{stab1} (cf. \cite[Theorem 2]{PV2006}) Assume that $\om = d_1 \tau$.  Then $E^{p,q}_{l,2} = 0$, if $1\le p \le q \le n-1$. 
If $0\le q \le n$, then $E^{0,q}_{l,2} = H^q_l (M^{2n})$. If $ 0\le p \le n$ then $E^{p,n}_{l+p,2} = H^{n+p}_{l+p} (M^{2n})$. Thus the  spectral sequence $\{ E^{p,q}_{l,r}, d_{l,r}\}$ stabilizes  at the term $E_{l,2}$. 
\end{theorem}

\begin{proof}  Assume that $\om= d_1 \tau$. Then  for any $d_{l-1}$-closed form $\alpha$ we have
\begin{equation}
d_1\tau \wedge \alpha = d_{l}(\tau \wedge \alpha). \label{exact0}
\end{equation}
Hence the induced operator  $\bar L$ in the exact sequence  (\ref{ex11}) satisfies
\begin{equation}
\bar L( H_{l-1} ^{q- (p+1)} (M^{2n})) = 0 \in  H_l ^{q+1 -p} (M^{2n}).
\label{exact0a} 
\end{equation}
 The equality (\ref{exact0a}) and the exact sequence  (\ref{ex11}) lead to the following exact  sequence for $0\le p \le q \le n-1$.
\begin{equation} 
0 \to  H^{q-p} _l(M^{2n}) \stackrel{\bar L^p}{\to} E^{p,q}_{l+p,1} \stackrel{\delta_{p,q}}{\to} H^{q-(p+1)} _{l-1}(M^{2n}) \to 0.\label{exact1}
\end{equation}

Using the isomorphism $E^{p,q}_{l+p, 1}  = H^{q-p}(\Pp^* (M^{2n}), d_l ^+)$ and the formulae (\ref{delta1}) and (\ref{lpex}) describing  $\delta_{p,q}$ and  $\bar L^p$,  we rewrite the
exact  sequence (\ref{exact1}) as follows
\begin{equation} 
0 \to  H^{q-p} _l(M^{2n}) \stackrel{[\Pi_{pr}]}{\to} H^{q-p}(\Pp^*(M^{2n}), d_l ^+) \stackrel{[d^-_l]}{\to} H^{q-(p+1)} _{l-1}(M^{2n}) \to 0.\label{exact1a}
\end{equation}

The proof of  the first and  second assertion of Theorem \ref{stab1}  is based on   our analysis  of the long exact sequence of cohomology groups
associated with  the short exact sequence (\ref{exact1a}).  By (\ref{stab1a}),  for $ 0 \le p \le q  \le n-1$, the differential  $d_{l+p,1} : E^{p, q}_{l+p, 1} \to E^{p+1, q}_{l+p,1}$ induces  the following boundary
operator
\begin{equation}
\hat d  ^-_l : H^{q-p}(\Pp^*(M^{2n}), d_l^+) \to H ^{q-p-1} (\Pp^*(M^{2n}), d_{l-1} ^+), \,  [\tilde \alpha]\mapsto  [d^-_l \tilde \alpha], \label{exact1a1}
\end{equation}
for $\tilde \alpha \in \Pp ^{q-p} (M^{2n})$.

\begin{lemma}\label{id1a}  The short exact  sequence (\ref{exact1a}) generates a short exact  sequence of the following chain complexes for $1\le p \le q \le n-1$:
\begin{eqnarray}
0 \to (H^{q-p}_l (M^{2n}), \tilde d_l : = 0) \stackrel{[\Pi_{pr}]}{\to} (H^{q-p}(\Pp^*(M^{2n}), d_l^+), \hat d_l^-)\to \nonumber\\
 \stackrel{[d^-_l]}{\to} (H^{q-(p+1)} _{l-1} (M^{2n}), \tilde d_{l-1} : = 0) \to 0.\label{id1a1}
\end{eqnarray}
\end{lemma}

\begin{proof}   It is  useful to  rewrite  the  sequence (\ref{id1a1})  of chain complexes as the following    diagram
$$
\xymatrix{
H^{q-(p-1)}_{l+1}(M^{2n})\ar[r] ^ {[\Pi_{pr}]} \ar[d] ^{\tilde d_{l+1} =0} &  H^{q-(p-1)}(\Pp^*(M^{2n}), d_{l+1} ^+) \ar[r]  ^{[d_{l+1}^-]} \ar[d]^{\hat d_{l+1} ^-}   & H^{q-p}_{l}(M^{2n})\ar[r] \ar[d] ^{\tilde d_l =0}  &\\
H^{q-p}_{l}(M^{2n})\ar[r] ^ {[\Pi_{pr}]} \ar[d] ^{\tilde d_l =0} &  H^{q-p}(\Pp^*(M^{2n}), d_{l} ^+) \ar[r]  ^{[d_{l}^-]} \ar[d]^{\hat d_l ^-}   & H^{q-(p+1)}_{l-1} (M^{2n})\ar[r] \ar[d] ^{\tilde d_{l-1} =0 } &\\
H^{q-(p+1)}_{l-1}(M^{2n})\ar[r] ^ {[\Pi_{pr}]} &  H^{q-(p+1)}(\Pp^*(M^{2n}), d_{l-1} ^+) \ar[r]  ^{[d_{l-1}^-]}  & H^{q-(p+2)}_{l-2}(M^{2n}) \ar[r] &
}
$$
To prove  Lemma \ref{id1a}, it suffices to show that the above diagram is commutative, or equivalently 
\begin{eqnarray}
\, \hat d^-_l[\Pi_{pr}] = \tilde d_l [\Pi_{pr}] = 0, \label{newcoh1a}\\
\, [d^-_{l-1}]\hat d^-_l  = \tilde d_{l-1} [d^-_l] = 0.\label{newcoh2a}
\end{eqnarray}

 Let $\alpha \in H^{q-p}_l(M^{2n})$ and $\tilde \alpha \in \Om^{q-p}(M^{2n})$  its representative.    Let $\tilde \alpha = \tilde \alpha _{pr} + L \tilde \beta _{pr} + L^2 \gamma$  be the Lefschetz decomposition of $\tilde \alpha$.  Using  $d_l\tilde \alpha = d^+_l\tilde \alpha = d^+_{l} \tilde \alpha _{pr} =0$ we obtain    $L d^- _l \tilde \alpha _{pr} = L d^+_{l-1}\tilde \beta _{pr} + L^ 2 (d^-_{l-1} \tilde \beta _{pr} + d_{l-2} \gamma)$. Hence
$$
0=L(d_l^-\tilde{\alpha}_{pr}+d_{l-1}^+\tilde{\beta}_{pr})+L^2(d_{l-1}^-\tilde{\beta}_{pr}+
d_{l-2}\gamma),
$$
which implies  $d_l^-\tilde{\alpha}_{pr}+d_{l-1}^+\tilde{\beta}_{pr}=0$ thanks to the uniqueness of the second Lefschetz  decomposition. It follows
\begin{equation}
\hat{d}_l^-[\Pi_{pr}]\alpha=\hat{d}_l^-[\tilde{\alpha}_{pr}]=[d_l^-\tilde{\alpha}_{pr}]=-[d_{l-1}^+\tilde{\beta}_{pr}]
=0\in H^{q-(p+1)}(\Pp^*(M^{2n}), d^+_{l-1}). \label{stab63a}
\end{equation}

Let $\beta \in H^{q-p}(\Pp^*(M^{2n}), d_l^+) $  and $\tilde \beta \in \Pp ^{q-p}(M^{2n})$ its representative. Then
\begin{equation}
\, [ d_{l-1}^- ] \hat d^-_{l} \beta = [d_{l-1} ^- d^-_{l} \tilde \beta] = 0 \in H^{q-(p+2)}_{l-2}(M^{2n}). \label{stab64a}
\end{equation}
Clearly   (\ref{newcoh1}) and (\ref{newcoh2})  follow from (\ref{stab63}) and (\ref{stab64a}).  This  completes the proof of Lemma \ref{id1a}.
\end{proof}

 The  short exact sequence  (\ref{id1a1}) in Lemma \ref{id1a} generates the following associated  long exact sequence of the cohomology groups
 \begin{equation}
 \to E^{p-1, q}_{l+p, 2} \to H^{q-p}_{l}(M^{2n}) \stackrel{\delta} \to H^{q-p}_{l}(M^{2n}) \to E^{p, q}_{l+p, 2} \to  H^{q-(p+1)}_{l-1}(M^{2n})  \stackrel{\delta} {\to}   \label{stab11}
 \end{equation}
where  $\delta$ is the connecting  homomorphism.

\begin{lemma}\label{deltata} We have $\delta (x) =  x$ for all $x \in H^{q-(p+1)}_{l-1}(M^{2n})$  and for all $1\le p \le n-1$.
\end{lemma}
\begin{proof}[Proof of Lemma \ref{deltat}]  Let $x  \in H^{q-(p+1)}_{l-1}(M^{2n})$.  Using (\ref{id1a1}) we write
$x= [d^-_l] y, y \in  H^{q-p}(\Pp ^*( M^{2n}),d^+_l )$. By definition of  the connecting  homomorphism
we have $\delta x = [\hat  d^-_l y] =  x$.  This completes the proof of Lemma  \ref{deltata}. 
\end{proof}

Clearly  Lemma \ref{deltata} implies the first assertion  of Theorem \ref{stab1}.

Now let us   consider   the case $ p = 0,\,   q\le n-1$. Then
$ E^{-1, q}_{l,1} = 0$.  
The previous  short exact sequence  (\ref{id1a1}) of chain complexes is now replaced by the new
one, where   the cohomology  groups    on the line  containing $E^{-1, q}_{l,1}$  left and right to $E^{-1, q}_{l,1}$  are  zero.  Let us write the new short exact sequence explicitly
as the following commutative diagram
$$
\xymatrix{
0\ar[r] ^ {0} \ar[d] ^{0} &  0 \ar[r]  ^{0} \ar[d]^{0}   & 0 \ar[r]  \ar[d] ^{0}  &\\
H^{q}_{l}(M^{2n})\ar[r] ^ {[\Pi_{pr}]} \ar[d] ^{\tilde d_l =0} &  H^{q}(\Pp^*(M^{2n}), d_{l} ^+) \ar[r]  ^{[d_{l}^-]} \ar[d]^{\hat d_l ^-}   & H^{q-1}_{l-1} (M^{2n})\ar[r] \ar[d] ^{\tilde d_{l-1} =0 } &\\
H^{q-1}_{l-1}(M^{2n})\ar[r] ^ {[\Pi_{pr}]} &  H^{q-1}(\Pp^*(M^{2n}), d_{l-1} ^+) \ar[r]  ^{[d_{l-1}^-]}  & H^{q-2}_{l-2}(M^{2n}) \ar[r] &
}
$$
The  resulting exact sequence of the cohomology  groups now are
$$ 0 \to H^q _{l}(M^{2n}) \to   E^{0, q}_{l, 2} \to  H  ^{q-1}_{l-1}(M^{2n}) \stackrel{ \delta} {\to } H ^ {q- 1} _{l -1}(M^{2n})...$$
Since $ \delta  = Id$, we obtain 
\begin{equation}
E^{0, q}_{l,2} = H^{q}_l (M^{2n}), \label{1411b}
\end{equation}
which proves the second assertion of  Theorem \ref{stab1}.

Next we compute $E^{p, n}_{l, 2} = E^{p, n}_{l, 1} / d_{l,1} (E^{p-1, n}_{l,1})$ for $0\le p \le n-1$.  Since $\om = d_1 \tau$,   the map
$[L] : H^{n -(p +2)}_{l-1} (M^{2n}) \to  C^{n-p}_l$  sends  $[\alpha]$ to $[d_l (\tau \wedge \tilde\alpha)] = 0 \in C^{n-p}_l$.  Thus Proposition \ref{prol1} implies that the following sequence is exact for $0\le p \le n-1$
\begin{equation}
0 \to C^{n-p}_l \stackrel{[L^p]}{\to}E^{p,n}_{l+p, 1} \stackrel{\delta_{p,n}}{\to} C^{n-(p+1)}_{l-1} \stackrel{[ L^{p+1}]}{\to} H^{n+p+1}_{l+p} (M^{2n})\to 0.\label{stab5}
\end{equation}
Set  for $ -1\le p \le n-1$
\begin{equation}
T^{n-(p+1)}_{l-1} : = \ker [L^{p+1}]:  C^{n-(p+1)} _{l-1} \to H^{n+p+1}_{l+p} (M^{2n}).\label{deft}
\end{equation}
 Then we obtain from the exact sequence  (\ref{stab5})  the following  short exact sequence
\begin{equation}
0 \to C^{n-p}_l \stackrel{[L^p]}{\to}E^{p,n}_{l+p, 1} \stackrel{\delta_{p,n}}{\to} T^{n-(p+1)}_{l-1}\to 0.\label{stab6}
\end{equation}

Using the isomorphism $E^{p,n}_{l+p, 1}  = \Pp ^{n-p} (M^{2n}) /  d^+_l (\Pp^{n-p-1}(M^{2n}))$ and the formulas  (\ref{new1a}) and (\ref{deltac}) describing $[L^p]$ and $\delta_{p,n}$
in the exact sequence  (\ref{prol1a}) of Proposition \ref{prol1}, we rewrite the short exact sequence (\ref{stab6})  as follows
\begin{equation}
0 \to C^{n-p}_l \stackrel{ [\Pi_{pr}]}{\to} \frac{\Pp ^{n-p}(M^{2n})}{  d^+_l (\Pp^{n-p-1}(M^{2n}))} \stackrel{[d_l ^-]}{\to} T^{n-(p+1)}_{l-1} \to 0. \label{stab61}
\end{equation}

Recall that  the map $[\Pi_{pr}]$  is already defined  in section 4.  It is  the quotient map of the map $\Pi_{pr}: (\ker d_l ^- \cap \Om ^{n-p} (M^{2n})) \to \Pp ^{n-p} (M^{2n})$, see  the explanation of (\ref{new1a}).

Next recall that the map $[d_l^- ]$   is the quotient map  of the  map $d_l ^-: \Pp ^{n-p}(M^{2n}) \to \Pp ^{n-p-1}(M^{2n}) \cap \ker d^-_{l-1}$, see the explanation of  (\ref{deltac}). (We now explain why this map is also  well-defined in (\ref{stab61}). First we have $d^-_l (d ^+ _l \gamma)  = d^+_{l-1} d^-_l \gamma =  d_{l-1} (d^-_l \gamma) = 0 \in C^{n-(p+1)}_{l-1}$. Furthermore for $\alpha \in \Pp ^{n -p} (M^{2n})$
$$ L^{p+1} d_l^- (\alpha) = L ^p ( Ld_l ^- \alpha ) = L^pd_l  \alpha = d_{l +p} L^p \alpha   = 0 \in H^{n +p +1} _{ l+p } (M ^{2n}).$$
Hence 
$$[d_l^-] (\frac{\Pp ^{n-p}(M^{2n})}{  d^+_l (\Pp^{n-p-1}(M^{2n}))})\subset  T ^{n -(p +1)}_{l-1} =  \ker  [L^{p+1}].$$ 
Thus $[d_l^-]$ is well-defined.)

Note that the differential $d_{l+p, 1} : E^{p,n}_{l+p, 1} \to E^{p+1, n}_{l+p, 1}$ induces the following boundary operator
\begin{equation}
\hat d^-_l: \frac{\Pp ^{n-p}(M^{2n})}{  d^+ _l(\Pp^{n-p-1}(M^{2n}))} \to \frac{\Pp ^{n-p-1}( M^{2n})}{  d^+_{l-1} (\Pp^{n-p-2}(M^{2n}))}, [\tilde \alpha] \mapsto [d^-_l \tilde \alpha], \label{stab62}
\end{equation}
for  $\tilde \alpha \in \Pp^{n-p}(M^{2n})$. The map $\hat d^-_l$ is well-defined, since by (\ref{com1a}) $d^-_l d^+_l \alpha  = d^+_{l-1} d ^- _l \alpha$ for $\alpha \in \Pp^{n-p-1}(M^{2n})$.

\begin{lemma}\label{newcoh} The short exact sequence (\ref{stab61})  generates a short exact sequence of the following chain complexes 
\begin{equation}
\xymatrix{
0 \ar[r] & C^{n-p}_l \ar[d]^{\tilde d_{l}:= 0} \ar[r]^{ [\Pi_{pr}]} & \frac{\Pp ^{n-p} (M^{2n})}{  d^+_l (\Pp^{n-(p+1)}(M^{2n}))}\ar[d]^{ \hat d^-_l} \ar[r]^{[d_l ^-]}& T^{n-(p+1)}_{l-1}\ar[d]^{ \bar d_{l-1}: = 0}\ar[r] &  0\\
0 \ar[r] & C^{n-(p+1)}_l  \ar[r]^{ [\Pi_{pr}]} & \frac{\Pp ^{n-(p+1)} (M^{2n})}{  d^+_l (\Pp^{n-(p+2)}(M^{2n}))}\ar[r]^{[d_l ^-]} & T^{n-(p+2)}_{l-1}\ar[r]&  0
}
\label{stab6a}
\end{equation}
\end{lemma}
\begin{proof}[Proof of Lemma \ref{newcoh}]  It suffices to show that 
\begin{eqnarray}
\, \hat d^-_l[\Pi_{pr}] = \tilde d_l \Pi_{pr} = 0, \label{newcoh1}\\
\, [d^-_{l-1}]\hat d^-_l  = \bar d_{l-1} [d^-_l] = 0.\label{newcoh2}
\end{eqnarray}

 Let $\alpha \in C^{n-p}_l$ and $\tilde \alpha \in \Om^{n-p}(M^{2n})$  its representative.    Let $\tilde \alpha = \tilde \alpha _{pr} + L \tilde \beta _{pr} + L^2 \gamma$  be the Lefschetz decomposition of $\tilde \alpha$.  Using  $d_l\tilde \alpha = d^+_l\tilde \alpha$ we obtain    $L d^- _l \tilde \alpha _{pr} = L d^+_{l-1}\tilde \beta _{pr} + L^ 2 (d^-_{l-1} \tilde \beta _{pr} + d_{l-2} \gamma)$. Hence
 \begin{equation}
\hat d^-_l [ \Pi_{pr}]\alpha=  [d^-_l \tilde \alpha _{pr}]=  [d_{l-1}^+ \tilde \beta _{pr}] = 0 \in  \frac{\Pp ^{n-p-1}( M^{2n})}{  d^+_{l-1} (\Pp^{n-p-2}(M^{2n}))}.\label{stab63}
 \end{equation}
Let $\beta \in (\Pp ^{n-p}(M^{2n}))/ d^+_l (\Pp^{n-p-1}(M^{2n}))$  and $\tilde \beta \in \Pp ^{n-p}(M^{2n})$ its representative. Then
\begin{equation}
\, [ d_{l-1}^- ] \hat d^-_{l} [\beta] = [d_{l-1} ^- d^-_{l} \tilde \beta] = 0. \label{stab64}
\end{equation}
Clearly   (\ref{newcoh1}) and (\ref{newcoh2})  follow from (\ref{stab63}) and (\ref{stab64}).  This  completes the proof of Lemma \ref{newcoh}.
\end{proof}

 The  short exact sequence  (\ref{stab6a}) in Lemma \ref{newcoh} generates the following associated  long exact sequence of the cohomology groups
 \begin{equation}
E^{p-1,n}_{l+p,2} \to  T^{n-p}_{l} \stackrel{\delta} \to C^{n-p}_{l} \to E^{p, n}_{l+p, 2} \to  T^{n-(p+1)}_{l-1}  \stackrel{\delta} {\to}  C^{n-(p+1)}_{l-1}  \to , \label{stab11}
 \end{equation}
where  $\delta$ is the connecting  homomorphism.

\begin{lemma}\label{deltat} We have $\delta (x) =  x$ for all $x \in T^{n-(p+1)}_{l-1}$  and for all $0\le p \le n-1$.
\end{lemma}
\begin{proof}[Proof of Lemma \ref{deltat}]  Let $x  \in T^{n-(p+1)}_{l-1}$.  Using (\ref{stab6a}) we write
$x= [d^-_l] y, y \in  (\Pp ^{n-p}( M^{2n}))/ d^+_l (\Pp^{n-p-1}(M^{2n}))$. By definition of  the connecting  homomorphism
we have $\delta x = [\hat  d^-_l y] =  x$.  This completes the proof of Lemma  \ref{deltat}. 
\end{proof}

Now let us complete the proof of Theorem \ref{stab1}. From Lemma \ref{deltat} and  the long exact sequence (\ref{stab11}) we obtain $E^{p,n}_{l+p, 2} = C^{n-p}_{l}/ T ^{n-p}_{l}$. Taking into account  (\ref{deft}) which defines $T^{n-p}_{l}$ to be the kernel of the surjective homomorphism
$[L^{p}]: C^{n-p}_{l} \to  H^{n+p}_{l+p}(M^{2n})$, we    conclude that  
\begin{equation}
E^{p,n}_{l+p, 2} = H^{n+p}_{l+p}(M^{2n})\text { for } 0\le p \le n-1.\label{stab22n}
\end{equation}
 Next, by Lemma \ref{fk1} $E^{n,n}_{l, 1} = C^\infty (M^{2n})$.  Since $d_{l, 1} (E^{n,n}_{l,1}) = 0$, using (\ref{stab62}) we get
\begin{equation}
E^{n, n}_{l+n, 2} = \frac{ C^\infty (M^{2n})}{d_l ^- (\Pp ^1 (M^{2n}))} =  H ^0 (\Pp ^* (M^{2n}), d_l ^-).\label{stab22p}
\end{equation}
By Proposition  \ref{primi}, $ d_l ^- $ is proportional to $ (d_l)^*_\om$. Applying the symplectic star operator we get from (\ref{stab22p})
\begin{equation}
E^{n, n}_{l+n, 2}= H ^0 (\Pp ^* (M^{2n}), d_l ^-)= H^0 (\Pp^*(M^{2n}), (d_l) ^* _\om) =  H^{2n}_{n +l} (M^{2n}).\label{stab22m}
\end{equation}
Clearly
the third assertion of Theorem \ref{stab1} follows from (\ref{stab22n}) and (\ref{stab22m}).  The last assertion of Theorem \ref{stab1}  follows immediately.
This completes the proof of Theorem \ref{stab1}.  
\end{proof}

From the exact sequence (\ref{exact1})  we  obtain immediately  the following

\begin{corollary}  Assume  that $\om = d_1 \tau$. For $0\le p \le q\le n-1$ we have 
\begin{equation}
E^{p, q}_{l+p, 1} = H^{q-p}_l (M^{2n}) \oplus H^{q-p-1}_{l-1} (M^{2n}).\label{stab12}
\end{equation}
\end{corollary}

Theorem \ref{stab1}  can be  generalized as follows.   Assume that $\om ^p = d_T \rho$  for some $\rho \in \Om^{2p-1} (M^{2n})$, in particular $d_T \om ^p = 0$. Clearly $d_p (\om ^p ) = 0 = d_T (\om ^p)$ implies that $ T = p$, since $L ^p : \Om^1 (M^{2n}) \to \Om ^{1+2p} (M^{2n})$ is injective.  Furthermore we have
$\om^{k+T} = d_{k +T}(\rho \wedge \om)$ for  all $ t\ge 0$. Thus there exists a minimal number $T$ such that
$\om ^T = d_T \rho$  for some $\rho \in \Om^{2T-1} (M^{2n})$.

\begin{theorem}\label{nil} (cf. \cite[Theorem 3]{PV2006}) Assume that $\om ^T = d_T \rho$  and $T \ge 2$. Then  the spectral sequence $(E^{p, q}_{l, r}, d_{l, r})$ stabilizes at the term $E^{*,*}_{l, T+1}$. 
\end{theorem}

\begin{proof} 
Our  proof of Theorem \ref{nil} exploits  the construction of   the exact couple   associated  with a filtered  complex. We use many ideas from \cite{Pietro2006}.  The main idea is to find a short  exact sequence, whose middle term is $E^{*,*}_{l, T}$, and moreover, this  short exact sequence is induced from  the trivial action of the  operator $L ^T$ on  (a part of)  complexes entering in the derived exact couples  (cf. with  the  proof of Theorem \ref{stab1}). The condition $T \ge 2$  is necessary for  Lemma \ref{isom1} below.

Let us begin with  recalling  the construction of  the derived  exact couple associated with a
filtered complex  $(F^pK^*, d_l)$ following \cite[p.37-43]{McCleary2001}.
We associate with  a  filtration $(F^pK^*, d_l)$   the 
following exact couple  
\begin{equation}\xymatrix{ D^{p+1, *}_l \ar[r]^i & D^{p, *}_l \ar[d]^j\\
 &  E^{p, *}_{l, 1}\ar[lu]^\delta },\label{excp}
 \end{equation}
where $D^{p, q}_l: = H^{p+q} (F^pK^*, d_l)$, which  we also abbreviate as $H^{p+q}_l (F^p K^*)$, and
$$\to D^{p+1, q-1}_l \stackrel{i}{\to }D^{p,q}_l\stackrel{j}{\to} E^{p, q}_{l,1}\stackrel{\delta} {\to} D^{p+1, q}_l\stackrel{i}{\to} D^{p, q +1}_{l,1} \stackrel{j}{\to}E^{p,q+1}_{l,1} \to$$
is the long exact  sequence of cohomology groups associated with the following short exact sequence of chain complexes  
\begin{equation}
0 \to (F^{p+1} K^{p+q}, d_l) \stackrel{\tilde i}{\to} (F^{p} K^{p+q}, d_l) \stackrel{\tilde j}{\to } (E^{p, q} _{l, 0}, d_{l,0}) \to 0.
\end{equation}
The differential $d_{l, 1}:E^{p, q}_{l, 1} \to E^{p+1, q}_{l, 1}$, defined in (\ref{stab1f}),  satisfies the following relation: $d_{l, 1} = j\circ \delta$.
 We refer the reader to \cite{McCleary2001} for  a comprehensive exposition  on the relation between the spectral sequence of a filtration and its associated exact couple.


 Set  $(D^{*,*} _{l})^0 : = D^{*,*} _{l}$. We define the   $t$-th derived  exact couple  of  the exact couple (\ref{excp}), $t\ge 1$,
\begin{equation}
(D_l^{p+1,q-1})^{(t)}\stackrel{i^{(t)}}{\to}(D_l^{p-t, q+t})^{(t)}\stackrel{j ^{(t)}}{\to} E^{p,q}_{l,t+1}\stackrel{\delta^{(t)}}{\to}(D_l^{p+1,q})^{(t)}\label{dct}
\end{equation}
 inductively as follows \cite[p. 38]{McCleary2001}.
\begin{eqnarray}
(D_l^{p,q}) ^{(t)}: = i  (D_l^{p+1, q-1})^{(t-1)} \subset D_{l} ^{p, q}, \label{excdt}\\
i ^{(t)} (i ^t x)  : = i  (i^t x) \text { for } i ^t x \in (D_l^{p, q}) ^{(t)}, \label {excit}\\
E^{p, q}_{l, t+1} : = \frac{\ker d_{l, t}  \cap E^{p, q}_{l, t}}{d_{l, t-1} (E^{p -t+1, q+t -2}_{l, t})}, \label{exet}\\
j ^{(t)} (i ^t x): = [j^{(t-1)} \circ (i^{(t-1)} x )], \label{excjt}\\
\delta ^{(t)} ([e]): = \delta ^{(t-1)} (e) \in i (D_l^{p, q})^{(t-1)},\\
d_{l, t+1} : = j ^{(t)} \circ \delta ^{(t)}.
\end{eqnarray}

Next we consider the following commutative diagram  
\begin{equation}
\xymatrix{D^{0,q-p}_l \ar[r]^{ \bar L}\ar[d]^{\bar L^p} &  D^{0, q-p+2}_{l+1}\ar[d] ^{\bar L^{p-1}}\\
D^{p, q}_{l+p} \ar[r] ^{i} & D^{p-1, q+1}_{l+p} 
} \label{excp1}
\end{equation}
where $\bar L ^p$ is induced by  the linear  operator $L^p: \Om ^{q-p} (M^{2n}) \to \Om ^{q+p} (M^{2n})$.

The diagram (\ref{excp1})  leads us to consider the following  diagram
\begin{equation}
\xymatrix{(D^{0,q-p}_l) ^{(t)} \ar[r]^{\bar L^{(t)}}\ar[d]^{(\bar L^p) } &  (D^{0, q-p+2}_{l+1})^{(t)}\ar[d] ^{(\bar L^{p-1}) }\\
(D^{p, q}_{l+p})^{(t)} \ar[r] ^{i^{(t)}} & (D^{p-1, q+1}_{l+p}) ^{(t)} 
}\label{excpk}
\end{equation}
where $\bar L^{(t)}$  (resp. $(\bar L ^p)$) is the restriction of $\bar L$  (resp. $\bar L^p$)  to $(D^{p, q}_{l+p})^{(t)}$.

\begin{lemma}\label{isom1} For $t \ge 1$ and $p\ge 1$, $T\ge 2$,  the following statements hold.\\
1. The diagram (\ref{excpk}) is commutative.\\
2. $(\bar L^{p})  $ is an isomorphism.\\
3. $\im \, (i ^{(t)})  = \im \,( (\bar L ^{p-1}))$. \\
4. $(\bar L^{(T-1)}) (D^{0,q-p}_{l}) ^{(T-1)} = 0$, if $d_T\om ^T = 0$.
\end{lemma}

\begin{proof} 1. The commutativity of (\ref{excpk}) is an immediate   consequence of  the commutativity of the diagram
(\ref{excp1}).

2. We  prove the second assertion of Lemma \ref{isom1} by induction, beginning with  $t = 1$.    Let $x = i (\alpha) \in  (D^{p, q} _{l+p}) ^{(1)} =  i (D^{p+1, q-1}_{l+p}) = i (H^{p+q}_{l+p}(F^{p+1}K^*)) \subset D^{p,q}_{l+p}$. Then  there is an element $\alpha' \in \Om ^{q-p-2} (M^{2n})$ such that   $[L^{p+1}\alpha'] = \alpha  \in D^{p+1, q-1}_{l+p}$, or equivalently $d_{l+p} (L^{p+1} \alpha') = 0$. Hence $L^{p+1} (d_{l-1} (\alpha')) = 0$. Since $d_{l-1}\alpha ' \in \Om ^{q-p-1}(M^{2n})$, $L^{p+1} (d_{l-1} (\alpha')) = 0$ implies  that $d_{l-1} \alpha ' = 0$, so  $\alpha' \in H^{q-p-2}_{l-1} (M^{2n})$. 
Hence $x = \bar L^p ( L \alpha')$,  and $L\alpha ' = i (\alpha') \in  (D^{0, q-p}_l )^{(1)}$.  This proves that  the linear map $(\bar L^{p})$
is  surjective  for $t=1$.  Furthermore, the map $(\bar L^{p})$  is  injective  for $t =1$, since  $L^p: \Om^{q-p}M^{2n} \to \Om ^{q+p}M^{2n}$  is  injective,  and $L^p d_{l} = d_{l+p} L^p$.  This proves Lemma \ref{isom1}.2  for $t =1$.

Now assume that Lemma \ref{isom1}.2  has been proved for $t = k$.  Since  $(D^{0, q-p}_l) ^{(k+1)}$ is a subset
of $(D^{0, q-p}_l) ^{(k)}$  the injectivity of $(\bar L^{p})$   follows  from the inductive statement. 
The surjectivity of $(\bar L^p)$ also  follows from the  commutativity of the diagram (\ref{excpk}), which implies
that  $(\bar L^{p-1})$ maps  the image $ i (( D^{0,q-p}_l)^{(k)})$  onto the set $i^{(k)} ( D^{p, q}_{l+p}) ^{(k)} 
= (D^{p-1, q+1}_{l+p})^{(k+1)}$.
This proves  Lemma \ref{isom1}.2 for  all  $t\ge 2$.

3. Clearly Lemma \ref{isom1}.3 is a consequence  of Lemma \ref{isom1}.2  and the commutativity of the diagram (\ref{excpk}).

4. Let us compute $\bar L ^{(T-1)}\beta$ for  $\beta \in (D^{0, q-p-2}_l)^{(T-1)}$, where $T\ge 2$.  By definition
$\bar L ^{(T-1)} \beta = \bar L (i^{(T-1)}[\tilde \beta])$, where  $\tilde \beta = L^{T-1}\hat \beta\in \Om^{q-p-2}(M^{2n}) $ for some $\hat \beta \in \Om^{q-p - 2T}(M^{2n})$, and $[\tilde \beta] \in D^{0,q-p}_l$, in particular we have $d_l \tilde \beta = 0$. Thus  $L^{T-1}d_{l-T-1} \tilde \beta = 0$. Since $L^{T-1}: \Om^{q-p-2T +2 } (M^{2n}) \to \Om^{q-p}(M^{2n})$ is injective we get $d_{l-T +1} \hat \beta = 0$.  Now we have
\begin{equation}
\bar L^{(T-1)} \beta = \bar L ( i ^{(T-1)}\tilde \beta) = i^{(T-1)}([ L^T \hat \beta])\in  (D^{p-1, q+1})^{(T)}_{l+p}, \label{isom14}
\end{equation}
by Lemma \ref{isom1}.1.
 
Note that 
\begin{equation}
[L^T \hat \beta] = [d_{T} \rho \wedge \hat \beta] = [d_{T +(l-T+1)} ( \rho \wedge \hat \beta)] = 0\in D^{0, q-p +2}_{l+1}. \label{isom14a}
\end{equation}
 Clearly (\ref{isom14}) and (\ref{isom14a}) imply  the last assertion of Lemma \ref{isom1}.
\end{proof}

Lemma \ref{isom1}.4 and the $(T-1)$th-derived exact couple yield the following short exact sequence 
\begin{equation}
0 \to (D^{p-(T-1), q+(T-1)}_l) ^{(T-1)}  \stackrel{j^{(T-1)}}{\to}  E^{p, q}_{l,T} \stackrel{\delta^{(T-1)}}{\to} (D^{p+1, q}_l) ^{(T-1)} \to 0.\label{vani2}
\end{equation}

\begin{lemma}\label{complex}  The short exact sequence (\ref{vani2}) generates    a short exact sequence of the following chain complexes
$$ 
\xymatrix{
0 \ar[r] & (D^{p-(T-1), q+(T-1)}_l) ^{(T-1)}  \ar[r]_{j^{(T-1)}}\ar[d]^{\tilde d_{l, T} = 0} &  E^{p, q}_{l,T} \ar[r]_{\delta^{(T-1)}}\ar[d] ^{d_{l, T}} & (D^{p+1, q}_l) ^{(T-1)} \ar[d]^{\tilde d_{l,T} = 0}\ar[r] & 0 \\
0 \ar[r] & (D^{p+1, q}_l) ^{(T-1)}  \ar[r]_{j^{(T-1)}} & E^{p+ T, q-T +1}_{l,T} \ar[r]_{\delta^{(T-1)}}& (D^{p+T+1, q-T +1}_l) ^{(T-1)} \ar[r] & 0 
}
$$
\end{lemma}

\begin{proof}  It suffices to show that
\begin{eqnarray}
d_{l, T} j^{(T-1)} = 0 , \label{complex1}\\
\delta^{(T-1)} d_{l, T} = 0.\label{complex2}
\end{eqnarray}

The equality $d_{l, T} j^{(T-1)} = 0$ holds, since $d_{l, T}$ is a  quotient map of the linear operator $d_l$  acting
on $\Om^*(M^{2n})$, and $j^{(T-1)}$ associates a cycle in $D^{p-(T-1), q+(T-1)}_l\subset H^{p+q}_l (M^{2n})$ to its  class
in $E^{p,q}_{l, T}$.

The equality $\delta^{(T-1)} d_{l, T} = 0$ holds,   since $\delta ^{(T-1)}d_{l, T} =\delta ^{(T-1)} j^{(T-1)} \delta ^{(T-1)} = 0$.  This completes the proof of  Lemma \ref{complex}.
\end{proof}

Let us continue the proof of Theorem \ref{nil}.  From Lemma \ref{complex} we obtain the following long exact sequence 
of the associated cohomology groups
\begin{equation}
(D^{p-(T-1), q+(T-1)}_l)^{(T-1)}\stackrel{j^*}{\to} E^{p,q}_{l, T+1 }\stackrel{\delta^*}{\to} (D^{p+1, q}_l) ^{(T-1)} \stackrel{\p}{\to}
(D^{p+1, q}_l) ^{(T-1)} \label{vani3}
\end{equation}

\begin{lemma}\label{isom2}  For  $0\le p\le q\le n$  the connecting homomorphism  $\p:  (D^{p+1, q}_l) ^{(T-1)}\to (D^{p+1, q}_l) ^{(T-1)}$ in (\ref{vani3}) is equal to the identity.
\end{lemma}

\begin{proof} By (\ref{vani2})   for any  $x\in (D^{p+1, q}_l) ^{(T-1)}$  there exists $e\in E^{p,q}_{l,T}$ such that $y = \delta^{(T-1)} (e)$. Since $d^{p, q}_{l,T} (e) \in \ker \delta ^{(T-1)}$ there exists  an element $y \in (D^{p+1, q}_l)^{(T-1)}$
such that $j^{(T-1)} (y) =  d^{p,q}_{l, T}( e) = j ^{(T-1)} \delta^{(T-1)} (e)$.  Since $j ^{(T-1)}$ is injective, $y = \delta ^{(T-1)} (e)$.  By definition $\p (x) = y = \delta ^{(T-1)} (e)$.  It follows
that $\p (\delta ^{(T-1)} e) = \delta ^{(T-1)} e$. This completes the proof of Lemma \ref{isom2}.
\end{proof}

\begin{corollary} \label{nilt} For $ p\ge T$ we have $E^{p,q}_{l, T +1} = 0$.
\end{corollary}
\begin{proof} For $p \ge  T$ Lemma 5.11 yields the following  exact sequence
$$(D^{*, *}_l)^{(T-1)}\stackrel{Id}{\to}(D^{*, *}_l)^{(T-1)}\stackrel{j^*}{\to} E^{p,q}_{l, T+1 }\stackrel{\delta^*}{\to} (D^{*, *}_l) ^{(T-1)} \stackrel{Id}{\to}
(D^{*, *}_l) ^{(T-1)}$$
which implies Corollary \ref{nilt} immediately.
\end{proof}

It follows from Corollary \ref{nilt} that $d_{l, T+1} : E^{p, q} _{l,  T+1} \to E^{ p +  T +1, q-T}_{l , T+1} =0$ for all  $p \ge 0$. This  completes the proof of  Theorem \ref{nil}.
\end{proof}

We end this section with   presenting   a proof of the following  stabilization theorem. 

\begin{theorem}\label{stcom} (cf. \cite[Theorem 4]{PV2006})  Assume that  $(M^{2n}, \om,\theta)$  is  a compact connected globally
conformal symplectic  manifold. Then the spectral sequence $(E^{p,q}_{l,k},d _{l,k})$ stabilizes at the $E^{*,*}
_{l,2}$-term.  
\end{theorem}

\begin{proof}  By Theorem \ref{conf}  it suffices to prove  Theorem \ref{stcom} for the case of a symplectic manifold $(M^{2n}, \om)$, i.e. $\theta = 0$. 
The proof we   present  here  uses many ideas in    the proof  of  Theorem 2 in Di Pietro's Ph.D. Thesis \cite{Pietro2006} stated for  connected compact symplectic manifolds.

By Lemma \ref{fk1} $E^{p, q}_{l, 1} = 0 = E^{p, q}_{l,k}$ if $q< p$  or $q>n$  for all  $k\ge 1$. Thus it suffices to examine the terms $E^{p,p}_{l, k}$, $E^{p, p+r}_{l, k}$,  for $0\le p\le n-r$, $r\ge 1$, $k\ge 2$.

\begin{lemma}\label{pp} Assume that ($M^{2n}, \om, \theta)$ is a compact l.c.s. manifold, and $0\le p \le n$.\\ 
1. If  $(M^{2n}, \om, \theta)$ is a globally conformal symplectic manifold, then $E^{p, p}_{l, k} = \R$ for
all $l$ and $k\ge  2$. Moreover
$E^{p,p}_{l,k}$ is generated by   the $p$-th power of  the symplectic form $\om$. \\
2. If $(M^{2n},\om, \theta)$ is  not  conformal equivalent to a symplectic manifold,
then $E^{p,p}_{l, k} = 0$  for all $ p\not = l$ and for all $k \ge 1$.
\end{lemma}

\begin{proof}[Proof of Lemma \ref{pp}] By (\ref{0pl})  if $0 \le p \le n-1$ then 
\begin{equation}
E^{p, p}_{l, 1} =  H^{0}_{l-p} (M^{2n}).\label{i00}
\end{equation}
By (\ref{e13a}) we obtain
\begin{equation}
E^{n,n}_{l, 1} = C^0_{l-n} = C^\infty (M^{2n}).\label{nn1}
\end{equation}
We get from (\ref{nn1}) and (\ref{stab1a})
\begin{equation}
E^{n,n}_{l, 2} = C^\infty (M^{2n})/ d_{l,1} (E^{n-1, n}_{l,1}) = C^\infty (M^{2n})/ d^-_{l-n} (\Pp^1(M^{2n})) = H^{0}(\Pp^*(M^{2n}), d^-_{l-n}), \label{nni}
\end{equation}
By Corollary \ref{lich0} 
\begin{equation}
 H_{0} (\Pp^*(M^{2n}), d^-_l) = H _{0}(\Pp^*(M^{2n}), (d_{l-n})^*_\om) =H^{2n}(\Om^*(M^{2n}), d_{l}).\label{nn2}
 \end{equation}
Note that $d_{l,k} (E^{n,n}_{l,k}) = 0$ and $\im \, d_{l,k} \cap E^{n,n}_{l,k} = 0$ for all $k \ge 2$.  Using (\ref{nni}) we get
\begin{equation}
E^{n,n}_{l,k}  = E^{n,n}_{l,2} =  H_{0} (\Pp^*(M^{2n}), d^-_l) \text{ for all }  k \ge 2. \label{nnil}
\end{equation}
Combining  (\ref{nnil}), (\ref{nn2}) with  Corollary \ref{lich0}  we obtain the  assertion  of Lemma \ref{pp}  for  the cases $ p =0$ or $ p = n$.

Now let us consider $E^{p,p}_{k, r}$ with $0 < p < n$.
By (\ref{0pl}) $E^{p, p}_{l,1} = H^0 _{l-p} (M^{2n})$.   

Let us first assume that $M^{2n}$ is globally conformal symplectic.  Using Theorem
\ref{conf} we  drop $l$ in the lower index   of $d_{l, r}$  and $E^{p, q}_{l, r}$. 
First we note that $E^{p,p}_0$ is generated by $\om^p$.  Since $E^{p,p}_k$ is a quotient of
$E^{p,p}_0$  and $[\om ^p] \in E^{p, p}_k$  for all $k\ge 1$, taking into account $[\om^n]\not = 0 \in E^{n,n}_k$, we  complete the proof of Lemma \ref{pp}.1.

Now let us assume that $[\theta] \not = 0\in H^1 (M^{2n})$. Then
Corollary \ref{lich0} asserts that $H^0_l (M^{2n}) = 0$ for all $l\not = p$.  It follows that  $H^{p,p}_{l, \infty} = 0$
for all $l\not = p$. This complete the proof of Lemma \ref{pp}.
\end{proof}

\begin{lemma}\label{pp1}  Assume that $(M^{2n}, \om)$ is a connected compact symplectic manifold. Then  for $0\le p \le n-2$ and $k \ge 2$ we have $E^{p, p+1}_{ k}  = H^1(M^{2n})$.
Furthermore $E^{n-1, n}_k = E^{n-1,n}_2$  for  all $k \ge 2$.
\end{lemma}

\begin{proof}[Proof of Lemma \ref{pp1}]  Using (\ref{stab1f})  and (\ref{stab1a}) we note that $d_1 : E^{0,1}_1 \to E^{1, 1}_1$ is equivalent to  the map $d^-: H^{1} (\Pp^*(M^{2n}), d^+) \to H^0 (M^{2n},\R) = \R$.  By (\ref{plus1b}) $H^{1} (\Pp^*(M^{2n}), d^+) = H^1 (M^{2n})$. Hence  
\begin{equation}
E^{0,1}_2 = H^1(M^{2n}).\label{e012}
\end{equation}
It follows that,  the image $d_k (E^{0,1}_k ) =0$ for all $k \ge 2$.  Thus we get from (\ref{e012})
\begin{equation}
E^{0,1}_k  = H^1 (M^{2n}) \text{ for all } k \ge 2.\label{01i}
\end{equation}
This proves  Lemma \ref{pp1} for $E^{0,1}_k$, $k \ge 2$.
Since the operator $L^p :\Om^ 1(M^{2n}) \to \Om ^{2p +1}(M^{2n})$ is injective for all $p\le n-1$, using Lemma \ref{pp} we get
\begin{equation}
E^{0,1}_k \cong E^{0,1}_k \wedge E^{p,p}_k \subset E^{p, p+1}_k \text{ for all } k \ge 2.\label{pp1a}
\end{equation}
Note that $E^{p,p+1}_k $ is a quotient of $E^{p, p+1}_2$, which is isomorphic to  $E^{0,1}_2$ by Proposition \ref{symm}. Taking into account (\ref{01i}) we obtain from (\ref{pp1a})
\begin{equation}
E^{p,p+1}_k \cong E^{0,1}_2 = H^1 (M^{2n}) \text{ for all } p \le n-2 \text { and } k \ge 2.\label{pp1b}
\end{equation}
This completes the proof the first assertion  of Lemma \ref{pp1}. The second assertion of Lemma \ref{pp1}  follows   from the  observation that $d_{2} (E^{n-1, n}_{2}) = 0 = \im\, d_{2} \cap E^{n-1,n}_{l,2}$, and $d_k (E^{n-1,n}_k) = 0 = \im \, d_k \cap  E^{n-1, n} _k$ for all $ k \ge 3$.
\end{proof}

\begin{lemma}\label{pp2}  Assume that $(M^{2n}, \om)$ is a connected compact symplectic manifold. Then $E^{n-2,n}_k = E^{n-2,n}_2$ for  all $ k \ge 2$. Furthermore, for $0 \le p \le n -3$   and $k \ge 2$ we have
\begin{equation}
E^{p, p+2}_k \cong   E^{0, 2}_2.\label{main1}
\end{equation}
\end{lemma}

\begin{proof}  First we note that  for  $k \ge 2$ 
\begin{eqnarray}
d_k (E^{n-2,n}_k) = 0, \nonumber\\
\im\, d_k  \cap E^{n-2,n}_k = 0.\nonumber
\end{eqnarray}
Hence
\begin{equation}
E^{n-2,n}_k = E^{n-2,n}_2  \text{  for  all } k \ge 2.
\end{equation}
This proves the first assertion of Lemma \ref{pp2}. 
Next we observe that $d_2 (E^{0,2}_2 ) = 0$ and $\im d_2 \cap E^{0,2}_ 2 = 0$. 
Hence  we get
\begin{equation}
E^{0,2}_k  = E^{0,2}_2 \text{ for all } k \ge 2.\label{02k}
\end{equation}
Now we assume that $0\le p \le n-3$. Since $L^{p}: \Om ^2 (M^{2n}) \to \Om ^{2 + 2p}(M^{2n})$ is injective, using Lemma \ref{pp} we get
\begin{equation}
E^{0,2}_k \cong E^{0,2}_k \wedge E^{p, p}_k  \subset E^{p, p+2}_k .\label{ppk}
\end{equation}
Since $E^{p,p+2}_k $ is a quotient group of $E^{p, p+2}_2$, which is isomorphic to $ E^{0,2}_2$ by Proposition \ref{symm}, using  (\ref{02k}) and (\ref{ppk}) we get
\begin{equation}
E^{p,p+2}_k \cong E^{0,2}_2 \text{ for  all } 0 \le p \le n -3.
\end{equation}
This completes the proof of Lemma \ref{pp2}.
\end{proof}

\begin{lemma}\label{pp3} We have $E^{p,p +r}_k = E^{0, r}_2 $ for all $k\ge 2$, $ p + r \le n-1$  and $r \ge 3$.
Furthermore $E^{n-r,n}_k = E^{n-r,n}_2$ for all $k \ge 2$ and $ r\ge 3$.
\end{lemma}

\begin{proof}  We  prove Lemma \ref{pp3}  inductively  on $r$   beginning with $r= 3$.  For each $r$  we  will consider $E^{p,p +r}_k $  with $k$ and $p$ increasing inductively.
  First we note that
\begin{equation}
d_2 (E^{0,3}_2 ) = 0 \in  E^{2,2}_2, \label{222}
\end{equation}
since $E^{2,2}_2 = E^{2,2}_k $ for all $ k  \ge 2$ by Lemma \ref{pp}.   From (\ref{222})  we obtain  easily
\begin{equation}
E^{0,3}_k =  E^{0,3}_2 \text { for all } k \ge 2.\label{222a}
\end{equation}
Now   using  the injectivity of the map $L^p : \Om ^3(M^{2n}) \to \Om ^{2p +3}(M^{2n})$ for $ p \le n -4$ and Lemma \ref{pp},  we get from (\ref{222a}) 
\begin{equation}
E^{0,3}_2 = E^{0,3}_k = E^{0,3}_k \wedge  E^{p, p}_k  \subset E^{p, p +3}_k .\label{222b}
\end{equation}
Since $E^{p, p+3}_k $ is a quotient group of $E^{p, p+3}_ 2  = E^{0,3}_2$,  (\ref{222b}) implies
\begin{equation}
E^{p, p+3}_k = E^{0, 3}_2 \text { for  all } 0\le p \le n-4, \, k \ge 2. \label{222c}
\end{equation}
This proves the  first assertion of Lemma \ref{pp3} for $r=3$.
The second assertion of Lemma \ref{pp3} for $r =3$ follows from the  identities $\im d_k \cap  E^{n-3,n}_k = 0 $ and  $d_k ( E^{n-3,n}_k) = 0 \in E^{n +k-3,n-k +1}_k$ if $k \ge 2$, which is a consequence of Lemma
\ref{pp} if $k =2$. 

Repeating this procedure we have  for  each $ n \ge r\ge 3$ the following sequences  of identities with  $k \ge 2$  and $ 0 \le p \le n-r$. First by induction on $r$ we get
\begin{equation}
d_2 (E^{0, r}_2) = 0 \in E^{2, r-1}_2,\label{fini1}
\end{equation}
since $E^{2, r-1}_2 = E^{2,r-1}_k$ for all $k \ge 2$ by the induction step. From (\ref{fini1}) we obtain immediately
\begin{equation}
  E^{0, r}_k  = E^{0, r}_2  \text { for  all }  k \ge 2\label{fini2}
  \end{equation}
Since the map $L^p : \Om ^r(M^{2n}) \to \Om ^{2p +r}(M^{2n})$ for $ p \le n -r$ is injective, using Lemma \ref{pp},  we get from (\ref{fini2}) 
\begin{equation}
E^{0,r}_2 = E^{0,r}_k = E^{0,r}_k \wedge  E^{p, p}_k  \subset E^{p, p +r}_k .\label{fini3}
\end{equation}
Since $E^{p, p+r}_k $ is a quotient group of $E^{p, p+r}_ 2 $ which is isomorphic to $ E^{0,r}_2$  if $p +r \le n-1$ by Proposition \ref{symm},  (\ref{fini3}) implies
\begin{equation}
E^{p, p+r}_k = E^{0, r}_2 \text { for  all } 0\le p \le n-r-1, \, k \ge 2. \label{fini4}
\end{equation}
Thus we get
\begin{equation}
 E^{0, r}_k = E^{p, p+r}_k \text { for  all } 0\le p \le n-r-1, \, k \ge 2. 
\end{equation}
This completes the proof of the first assertion of  Lemma \ref{pp3}.

The second assertion of Lemma \ref{pp3} for the inductive $r$ follows from the  identities $\im d_k \cap  E^{n-r,n}_k = 0 $ and  $d_k ( E^{n-3,n}_k) = 0 \in E^{n +k-r,n-k +1}_k$ if $k \ge 2$, which is a consequence of the induction assumption that $ E^{n+k-r, n-k+1}_k =  E^{n+k-r, n-k+1}_2$ for $0\le k \le r-1$. 
\end{proof}
Clearly   Theorem \ref{stcom} follows from Lemmata \ref{pp}, \ref{pp1}, \ref{pp2} , \ref{pp3}.
\end{proof}

\begin{remark}\label{tsey}  1.   Our
stabilization theorem \ref{stcom}  gives an answer  to the Tseng-Yau question on the relation between
the group $H^{p} (M^{2n}, d^+) = E^{0, p}_1 (M^{2n})$ for $ 0 \le p \le n-1$ and the cohomology groups $H^* (M^{2n}, \R)$.

2. In the next section we  show  that if $(M^{2n}, \om)$ is  a compact K\"ahler manifold, then
the spectral sequence   stabilizes already at $E_1$-terms, see Theorem \ref{kaehler1}.  
\end{remark}

\section{K\"ahler manifolds}

In this section we prove that
the spectral sequence  $E^{p,q}_r$ stabilizes  at the  term  $E^{*,*}_1$, if $(M^{2n},J, g)$ is a compact K\"ahler 
manifold and $\om$ is the associated symplectic form  (Theorem \ref{kaehler1}).

Let  $(M^{2n}, J, g, \om)$ be a  compact K\"ahler manifold.  As before, denote by $d^*$ the formal adjoint  of $d$. Since the operator
$L$ commutes with the Laplacian $\triangle  _d : = d d^* + d^* d$ we get the induced Lefschetz decomposition
of the space of harmonic  forms on $M^{2n}$, and hence the  induced Lefschetz  decomposition  of $H^*(M^{2n}, \R)$. Let us denote  by $P H^q(M^{2n}, \R)$ the subset of primitive cohomology classes in $H^q(M^{2n}, \R)$. Note that
each primitive cohomology class $[\alpha] \in P H^q (M^{2n})$ has a  representative  which is $\triangle_d$-harmonic and primitive.

\begin{proposition}\label{pmhk}(\cite[Proposition 3.18]{TY2009}) Assume that $(M^{2n}, J,g, \om)$ is  a compact  K\"ahler manifold. For $q \le n-1$ we have
$H_q (\Pp^*(M^{2n}), (d)^*_\om )  = P H^q (M^{2n}) = H^{q} (\Pp^*(M^{2n}), d^+ )$.
\end{proposition}

\begin{proof} We give here another proof using \cite{Bouche1990}. Bouche proved that  if $(M^{2n}, J, g, \om)$ is a compact K\"ahler manifold,  the  coeffective cohomology groups  $ H^{2n-q} (\Cc^* M^{2n}, d)$ is  isomorphic  to the subgroup $H^{2n-q}_\om ( M^{2n}) : = \{ x \in H^{2n-q }(M^{2n}, \R)|\, Lx = 0\}$ for $0\le q \le n-1$ \cite[Proposition 3.1]{Bouche1990}. Next, using \cite[Corollary 2.4.2]{Brylinski1988},  or  the following  formula: $\Jj d ^*  \Jj ^{-1} =  d ^* _\om$, which is proved in the similar way as (\ref{conj2}) replacing (\ref{formal10}) in the  proof with (\ref{formal4}), we observe that the symplectic star operator $*_\om$ maps $H ^{2n-q}_\om (M^{2n})$ isomorphically onto the group  $PH^q (M^{2n})$.
As we have noted 
in Remark \ref{his1}.1,  the coeffective cohomology group  $H^q(\Cc^*M^{2n}, d)$ is isomorphic to    the  primitive homology group $H_{q} (\Pp ^*M^{2n})$. On the other hand Tseng-Yau proved that  the  group $H(\Pp^*M^{2n}, (d)^-_\om)$ is  isomorphic to the group $H_{q}(\Pp^*(M^{2n}), d ^-)$  \cite[Lemma 2.7, part II]{TY2009} as well as to the group $H^{q,+} (M^{2n})$, \cite[Proposition 3.5, part II]{TY2009}, see also   Proposition \ref{primi} and Proposition \ref{dpm1} above.
Combining these observations we complete the proof of Proposition \ref{pmhk}.
\end{proof}

\begin{theorem}\label{kaehler1}  Assume that $(M^{2n}, J, g, \om)$ is  a compact K\"ahler manifold. Then the spectral
sequence  $E^{p, q}_r$ stabilizes at $E_1^{*,*}$.
\end{theorem}

\begin{proof}   Since $(M^{2n}, J,g, \om)$  is  a compact symplectic  manifold, by Theorem \ref{stcom}  the spectral sequence
$E^{p,q}_r$ stabilizes at  $E_2$-terms. Thus to  prove  Theorem \ref{kaehler1}  it suffices to show  that
all the differentials $d_1: E^{p, q}_1 \to  E^{p+1,q}_1$  vanish.  By (\ref{stab1a}), if $q\le n-1$ then $d_{ 1} : E^{p,q}_{ 1} \to E^{p+1, q}_{1}$ is defined by the image of $d^- \tilde \alpha$, $\tilde \alpha \in \Pp^{q-p} (M^{2n})$. In this case it suffices to show that
any element $[\alpha]\in H^{q-p} (\Pp^*(M^{2n}), d^+)$ has a  representative $\alpha \in \Pp ^{q-p}( M^{2n})$ such that
$d^- \alpha = 0$. By the Hodge theory  for $d ^+ $  there is  a harmonic representative $\alpha$ of $[\alpha]$  such that $(d^+) ^*  \alpha = 0$.    Lemma \ref{dpm} implies that for such  harmonic form $\alpha$ we have
$d ^- \alpha = 0$.  This  implies $ d_1 (E^{p,q}_1) = 0$ for $q\le n-1$.  It remains to consider
the image $d_1 (E^{p, n}_1)$. By (\ref{stab62}) it suffices to show that  any  element
$[\alpha] \in E^{p,n}_1$  has a representative $\alpha \in \Pp^{n-p}(M^{2n})$ such that $d^- (\alpha) =0$. Using the Hodge
theory for $d^+$ and (\ref{stab62})  we choose  $\alpha$ to be the harmonic form.  By Lemma \ref{dpm} $d  ^- (\alpha) =0$.  This completes the proof
of Theorem \ref{kaehler1}.
\end{proof}

\section{Examples}

In this section we   consider two  simple examples of compact l.c.s. manifolds and their  primitive cohomologies. The first example is a  nilmanifold of Heisenberg type \cite{Sawai2007b},  the second example is  a  4-dimensional  solvmanifold described in  \cite{ACFM1989}, \cite{Banyaga2007}, \cite{Sawai2007}, \cite{HK2011}.
We calculate    the primitive cohomology of  these examples (Propositions \ref{mil1},  \ref{solv1}).  We   study  some properties  of primitive cohomology groups of a l.c.s.  manifold,  which is a   mapping torus of  a co-orientation preserving contactomorphism (Proposition \ref{inv2}). We show that the 4-dimensional solvmanifold  is a mapping torus of  a  coorientation preserving contactomorphism of  a connected   contact 3-manifold, which is not  isotopic to the identity (Theorem \ref{new1}). 

\medskip

Let $H(n)$ denote the $(2n+1)$-dimensional Heisenberg  Lie group and $\Gamma$  its lattice. It is well-known that the nilmanifold $N^{2n+2} : = (H(n)/\Gamma) \times S^1 $ has a canonical l.c.s. form $\Om$,  which we now describe following \cite{Sawai2007b}.   Note that the  Lie algebra $\h (n) \oplus \R$ of $H(n) \times \R$
is given by $\la X_i, Y_i, Z, A : [X_i, Y_i] = Z\ra _\R$. We denote by $x_i, y_i, z, \alpha$ the  dual basis. 
Clearly $d\alpha = 0$ and $d\Om = \alpha \wedge \Om$. Here we use the same notations for the extension of   $X_i, Y_i, Z, A,  x_i, y_i, z, \alpha, \Om$ to the right-invariants vector fields or differential forms on  $H(n) \times \R$,
as well as for the descending  vector fields or differential forms on $N^{2n+2}$.

\begin{proposition}\label{mil1}  Let $(N^{2n+2}, \Om, \alpha)$ be the  l.c.s. nilmanifold  described above. All the Lichnerowicz-Novikov cohomology groups $H^*(\Om^*(N^{2n+2}), d_{k\alpha})$  vanish, if $k \not = 0$. Consequently for $k \not = 0$ all  the groups $E^{p,q}_{k, r}$, $r\ge 1$, of the associated spectral sequences vanish, unless $q = n$ and $r =1$. The group $E^{p,n}_{k,1}$ is infinite  dimensional
for all $0 \le p \le n$.
\end{proposition}

\begin{proof}  The first  assertion of Proposition  \ref{mil1} is  a consequence of a result due to Millionshchikov, who proved that the  Lichnerowicz-Novikov
 cohomology  groups $H^*  (\Om^*(M), d_\theta)$  of a compact nilmanifold $M$  always vanish unless $\theta$  presents a trivial cohomology class in $H^1 (M,\R)$ \cite[Corollary 4.2]{Millionshchikov2005}.
The second assertion of Proposition \ref{mil1} for $E^{p,q}_{k, 1}$ is a consequence of the first assertion, combining with Lemma \ref{fk1} and the exact sequence (\ref{ex11}).  Since $\Om = d_\alpha (z)$, applying Theorem \ref{stab1}  we obtain the second assertion from the first assertion combining with
the  particular case $ r =1$ proved above. The third assertion follows from Lemma \ref{fk1} and  from the ellipticity of the operators $d_k^+$, see the proof of  Proposition \ref{dpm1}. This completes the proof of Proposition \ref{mil1}.
\end{proof}

Now we  shall show an example of a l.c.s.  4-manifold $M_{k,n}$, which is an Inoue surface of type $S^-$,  whose primitive cohomologies are non-trivial,  and we will explain  an implication of this non-triviality.  The  4-manifold $M_{n,k}$ has been described in  \cite{ACFM1989}, \cite{Banyaga2007}, \cite{Sawai2007}, \cite{HK2011}. Here we follow  the exposition in  \cite{Banyaga2007}. 
Let $G_k$ be the group of matrices of  the form

$$
\left(\begin{array}{clrr}      e^{kz} & 0& 0& x \\  
0  & e^{-kz}  & 0 & y\\
0 & 0 & 1 &  z\\
0 & 0 & 0  & 1
\end{array}\right),
$$
where $x,y,z \in \R$  and $k \in \R$ such that $e^k + e^{-k} \in \Z \setminus \{2\}$. The group $G_k$ is a connected solvable Lie group with a basis of  right invariant 1-forms
\begin{equation}
 dx -kxdz, \: dy + kydz, \: dz. \label{1inv}
\end{equation}
 
There exists a discrete subgroup $\Gamma _k \subset G_k$ such that $N_k = G_k/\Gamma_k$ is compact. The basis (\ref{1inv}) descends to a basis of 1-forms $\alpha, \beta, \gamma$ on $N_k$. The forms $\gamma$ and $\alpha \wedge \beta$ are closed and their cohomology  classes generate $H^1 (M, \R)$ and $H^2 (M, \R)$ respectively.

Now let $\lambda \in \R$ be a number such that $\lambda [\alpha \wedge \beta] \in H^2 (M, \Z)$. For given $k,n$ denote by $M_{k,n}$ the total space of 
the $S^1$-principal bundle over $N_k$ with the Chern class $n\lambda [\alpha \wedge \beta]$.  Let $\eta$ be a connection form on $M_{k,n}$, equivalently
\begin{equation}
d\eta = n \lambda (\alpha \wedge \beta).\label{curv1}
\end{equation}
For simplicity we  will denote  the pull back to $M_{k,n}$ by the projection $M_{k,n} \to N_k$ of a form $\theta$ on $N_k$ again by $\theta$.
Banyaga  showed that $M_{k,n}$ possesses many interesting l.c.s. structures. Here we  consider only two l.c.s. forms $ d_{-k\gamma} \eta = n \lambda (\alpha \wedge \beta) - k \gamma \wedge \eta $  and $d_{k\gamma} \eta = n \lambda (\alpha \wedge \beta) + k \gamma \wedge \eta $  discovered by Banyaga \cite[Remark 2]{Banyaga2007}. Note that $M_{k,n}$ carries no symplectic structure, since $H^2 (M_{k,n}, \R) = 0$ \cite{ACFM1989}.  Since
$M_{k,n}$ is compact, the Hodge theory  applied to $d_{\pm \gamma}$ yields that $H^{2-i} (\Om ^*(M_{n,k}), d_{\pm k \gamma}) = H^{2+i}(\Om^*(M_{n,k}), d_{\pm k \gamma})$. Since $[\pm k\gamma] \not = 0 \in H^1(M_{n,k}, \R)$, the Lichnerowicz deformed
differential  $d_{\pm k\lambda}$ is not gauge equivalent to the  canonical differential $d$. Hence  $H^0(\Om ^*(M^{2n}), d_{\pm k \gamma})  = 0$.  Denote by $\Pp_{\pm}^*(M^{2n})$ the space of primitive forms corresponding to the
l.c.s. form $d_{\pm k\gamma} \eta$. Corollary \ref{lich0} yields that  $H^{0} (\Pp^*_{\pm}(M_{k,n}), d^+_{\pm lk\gamma})  = 0$ for all $ l \not = 0$, and $H^0(\Pp_{\pm} ^*(M^{2n}), d) = \R$.

\begin{proposition}\label{solv1} 1.  $H^1 (\Om ^*(M_{k,n}), d_{\pm k \gamma}) = \R$.\\
2. $H^2(\Om ^*(M_{k,n}), d_{\pm k \gamma} ) = \R$.\\
3. $  H^{1}   (\Pp^*_{\pm} (M_{k,n}), d^+_{\pm  k \gamma}) = \R^2$.
\end{proposition}

\begin{proof}  It is known that $M_{k,n}$  is a complete solvmanifold. Indeed,  the algebra $\g_{k,n}$  of the corresponding solvable group possesses the  basis  $(X, Y, Z, T)$ dual to $(\alpha, \beta, \gamma, \eta)$ with the following properties \cite{ACFM1989}, or see (\ref{solv2d}) and (\ref{solv3d}) below.
\begin{eqnarray} [X, Z] = k X, \, [X, Y] = -n \lambda T, \,  [Y, Z] = - k Y, \label{solv2}\\
\, [X, T] = [Y, T] = [Z, T] = 0.\label{solv3}
\end{eqnarray}
Using (\ref{solv2}) and  (\ref{solv3}) we observe that the Lie subalgebras  $\la T\ra_\R \subset \la T, X \ra _\R \subset  \la T, X, Y \ra _\R$ are ideals of  $\g_{k,n}$, so  $M_{k,n}$ is completely solvable.
Now we  apply  the result 
by Millionshchikov \cite[Corollary 4.1, Theorem 4.5]{Millionshchikov2005}, which reduces the computation of the Novikov
cohomology groups of a compact complete solvmanifold $G/\Gamma$ to the computation of the induced Novikov cohomology groups of the Lie algebra $\g$ of $G$.  For our computation it is useful to rewrite (\ref{solv2}) and (\ref{solv3}) in  the dual basis of $\g^*_{k,n}$, or using the explicit  formulae for $\alpha, \beta, \gamma, \eta$ given in (\ref{1inv}), (\ref{curv1}) above to obtain
\begin{eqnarray}
d\alpha =- k \alpha \wedge \gamma, \: d\beta = k \beta \wedge \gamma, \label{solv2d}\\
 d\gamma = 0, \: d\eta =  n \lambda (\alpha \wedge \beta). \label{solv3d}
\end{eqnarray}

1. Abbreviate  $d_{\pm k \gamma}$ as $d_{\pm k}$. Using (\ref{solv2d}) and (\ref{solv3d}) we get
\begin{eqnarray}
d_{\pm k} \alpha = (k \pm k) \gamma \wedge \alpha, \: d_{\pm k} \beta = (- k \pm k) \gamma \wedge \beta , \label{solv4d}\\
 \gamma = d_{\pm k } (\pm 1/k), \: d_{\pm} \eta =  n \lambda (\alpha \wedge \beta) \pm k\gamma \wedge \eta.\label{solv5d}
\end{eqnarray}
Using (\ref{solv4d}) and (\ref{solv5d}) it is easy to compute that $\alpha$ is a generator of $H^1  (\Om^*(M_{k,n}), d_{-k})$,  and $\beta$ is  a generator
of $H^1 (\Om^*(M_{k,n}), d_k)$. This proves the first assertion of Proposition \ref{solv1}.

2. For computing $H^2  (\Om ^*(M_{k,n}), d_{\pm k})$ we use (\ref{solv4d}), (\ref{solv5d}), and  the following formulae
\begin{eqnarray}
d_{\pm k} (\alpha \wedge \beta ) = \pm k \alpha \wedge \beta \wedge \gamma, \nonumber\\
  d_{-k} (\alpha \wedge \gamma)= 0, \: d_{k}  (\beta  \wedge \gamma ) = 0, \nonumber\\
d_{\pm k} (\alpha \wedge \eta) = (-k \pm k) \gamma \wedge \alpha \wedge \eta, \: d_{\pm k} (\beta \wedge \eta) = (k \pm k)\gamma \wedge \beta \wedge \eta, \nonumber \\ 
d_{\pm k} (\gamma \wedge \eta) = - n \lambda \alpha \wedge \beta \wedge \gamma.\nonumber
\end{eqnarray}

It is easy to see that  $\alpha \wedge \eta$   is a generator of $H^2 (\Om ^*(M_{k,n}), d_{-k})$ and $\beta \wedge \eta$ 
is a generator of  $H^2(\Om^*(M_{k,n}), d_k) $. This proves the  second assertions of  Proposition \ref{solv1}. 

3. The third  assertion of Proposition \ref{solv1}  is a consequence of the first assertion and Formula(\ref{plus1c}).
This completes the proof of Proposition \ref{solv1}.
\end{proof}

In the remaining part of this section we   study  some properties  of primitive cohomology groups of l.c.s.  manifolds associated with   a co-orientation preserving contactomorphism.  We show  that the  l.c.s. solvmanifold  studied  before  is an example of a  l.c.s. manifold associated with  a non-trivial  contactomorphism.

Let $(M^{2n+1}, \alpha)$ be a  co-orientable contact manifold and $f $ be a co-orientation preserving  contactomorphism of $(M^{2n+1},\alpha)$, i.e.  $f^* (\alpha) = e ^h \cdot \alpha$ for some
$h \in  C^\infty (M^{2n})$.
The mapping torus  $M_f^{2n+2} = (M\times [0,1])/( [x, 0] = [f (x) , 1])$ of  a contactomorphism $f$ is a fibration over $S^1$ whose fiber is $M^{2n+1}$.   Let us denote this  fibration by $\pi: M_f^{2n+2} \to S^1$ with
$\pi ^{-1} (s) = [M, s]$.
 Let $f_t: M_f^{2n +2} \to M_f^{2n +2}$ be a  1-parameter family of diffeomorphisms defined by:
$$f_t ([x, s]) = [x, s+t \mod 1] \text {  for } t \in \R.$$
In particular  $f_1  ([x, 0]) = [f(x), 0])$.  Let us also denote by $\alpha$  the  contact 1-form on $[M, 0]$ obtained by identifying $M$  with $[M, 0]$. 
Let $B$ be  the vector field on $M^{2n+2}_f$ defined by $B([x,s]) = (d/dt)_{|t= 0} f_t ([x, s])$.  Since  $f ^* (\alpha) =  e ^h \alpha$  the following   1-form $\tilde \alpha$  
\begin{equation}
\tilde \alpha (x, t)_{|\pi^{-1} (t)}: =  e ^{ - th (x)}f_{-t} ^ * (\alpha),  \, \tilde \alpha (B) = 0.\label{ext1}
\end{equation}
is well-defined on $M^{2n+2}_f$, moreover
$$ f _t ^* ( \tilde \alpha) _{| \pi ^{-1} (t)}  = \tilde \alpha _{|\pi ^{-1} (0)} \text  {  for all } 0\le  t \le 1. $$ 
Set $\theta : = \pi ^* (dt)$.

\begin{proposition}\label{mtc} (cf. \cite[Proposition 3.3.]{BK2010}) 1. Assume that $ (M^{2n +1}, \alpha)$ is a compact co-orientable  contact manifold  and $f$ is  a co-orientation preserving contactomorphism. There exists  a positive number $c_0$ such that $(M^{2n+2}_f, \om_c:=d\tilde \alpha + c\theta \wedge \tilde\alpha, c\theta)$ is a l.c.s. manifold for
all $ c\ge c_0$.  

2. Assume that $f$ preserves the contact 1-form $\alpha$. Then $(M^{2n+2}_f, \om:=d\tilde \alpha + \theta \wedge \tilde\alpha, \theta)$ is a l.c.s. manifold. 
\end{proposition}

\begin{proof} 1.  Clearly (\ref{ext1}) implies that $ rk \, d\tilde \alpha  \ge rk \, d\alpha = 2n$. Using this  we conclude that there exists a positive number $c_0$ such that   $rk\, d\om_c = 2n +2$ for all $ c\ge c_0$, since  $M^{2n +1}$ is compact.  Further, $d (c \theta  ) = 0$ and
$ \om_c = d _{c\theta} (\tilde \alpha)$. This proves that $(M^{2n +2} _f, \om _c,  c\theta)$ is a l.c.s. manifold. 

2.  Assume that $f ^* (\alpha )=   \alpha$.  Then  $\tilde \alpha ( [x, t])_{\pi ^{-1} (t)} = f ^* _{-t}   \alpha$. It follows that $rk\,  d\tilde \alpha = rk\, d\alpha = 2n$, and $ rk \, \om _1 = 2n +2$.
Hence $\om _1 = \om$ is a l.c.s. form, taking into account $\om = d_{\theta } \alpha$.
\end{proof}


\begin{proposition} \label{inv2} 1. Suppose that $f_0$ and $f_1$ are   co-orientation  preserving contactomorphisms  of a   compact co-orientable contact manifold $(M^{2n+1}, \alpha)$.    The l.c.s. manifolds  $M^{2n+2}_{f_0}$ and
$M^{2n +2}_{f_1}$ are diffeomorphic, if $f_0$ and $f_1$ are isotopic. For   sufficiently large  number $c$ the primitive cohomology groups of $(M^{2n +2}_{ f_0}, \om _c, c \theta)$ and  of $(M^{2n +2}_{ f_1}, \om _c', c \theta)$ are isomorphic.

2. Let  $\theta$ be the Lee form of the associated  l.c.s form on $M^{2n +1}_f$.  If $f$ is isotopic to the identity, the Lichnerowicz cohomology groups $H ^*(\Om ^*  (M^{2n+2}_f), d_{c \theta})$   are zero, for any $c \not = 0$. 
\end{proposition}

\begin{proof} The first  assertion of Proposition \ref{inv2}.1 is well-known.
The second assertion of Proposition \ref{inv2}.1 is a consequence of  Theorem \ref{conf}, observing that  $\om _c - \om '_c = d_{c\theta} (\tilde \alpha - \tilde \alpha ')$.

Finally   Proposition \ref{inv2}.2 follows  from the first assertion, combining with the fact that the l.c.s. manifold $(M^{2n+1} \times S^1, d_\theta \tilde\alpha, \theta)$ associated to the identity  mapping of the contact manifold $(M^{2n+1}, \alpha)$ has
 vanishing Lichnerowicz-Novikov groups, taking into account the K\"unneth formula and the  formula $H^* (\Om ^*(S^1), d_{c dt}) = 0$ if $c \not = 0$.
This completes the proof of Proposition \ref{inv2}.
\end{proof}

Now we shall show that our l.c.s. manifold $(M_{k,n}, d_{k\gamma}\eta, k\gamma)$ is  a mapping torus of a non-trivial  co-orientation preserving contactomorphism.
First we prove the following

\begin{proposition}\label{cri1} Assume that  $\gamma$ is a closed  1-form  on a compact smooth manifold  $M$.  If $[\gamma]\in H^1 (M, \Z)$ and $\gamma$ is now-where vanishing, then there
is a submersion $f : M \to S^1$ such that  $f^* (dt) = \gamma$, where $dt$ is the canonical  1-form on $S^1$.
\end{proposition}

\begin{proof}  We use Tischler's argument in \cite{Tischler1970}.
 Since $S^1$ is the Eilenberg-Maclane  space  there exists a map $f_1: M \to S^1$ such that  $f^*([dt]) =  [\gamma]$.  Without loss of generality we assume that
$f$ is a smooth map.  Hence we have $f^* (dt) =  \gamma + dh$ for some smooth function $h$ on $M$.  Now we observe that  $f_1 ^* (dt) + dh =  (f_1 + \Pi \circ h)  (dt)$, where $\Pi : \R \to S^1$ is the natural projection. Clearly
the map $f =  f_1 + \Pi\circ h$ is a submersion, since $\gamma$ is no-where vanishing. This completes the proof  of Proposition \ref{cri1}.
\end{proof}

Now we are ready to   show the following implication of  Proposition \ref{solv1}.

\begin{theorem}\label{new1}  The l.c.s. manifold $(M_{k,n}, d_{k\gamma}\eta, k\gamma)$   is  a mapping torus 
of  a coorientation preserving contactomorphism  $f$  of a 3-dimensional  connected contact manifold.  Moreover $f$  is not isotopic  to the identity.
\end{theorem}

\begin{proof} Since  $H^1 (M_{k,n}, \R) = \R$  \cite{ACFM1989}, and $d\gamma = 0$, there exists a positive number $p$ such that $p [\gamma]$ is a generator of   $H^1 (M_{k,n}, \Z) = \Z = Hom (H_1(M_{k,n}, \Z), \Z)$ \cite[Chapter VI, 7.22]{Dold1972}.  Applying Proposition \ref{cri1} we conclude that
$M_{k,n}$ is a fibration over $S^1$ whose fibers are the foliation $\Ff_1 : =\{\gamma = 0\}$, and $f ^* (dt) = p \cdot\gamma$, since $\gamma$ is nowhere vanishing.  Denote by $\pi: M_{k,n} \to S^1$ the  corresponding  fibration.
Note that the restriction of $\eta$ to each  fiber $\pi^{-1}(t)$, $ t \in S^1$,  is a contact form, since $X, Y, T$ are tangent to the fiber and we have $\eta (T) = 1$,  $d\eta (X, Y)\not = 0$. 

First we will show that the fiber $F : = \pi^{-1}(t)$, $t \in S^1$,  is connected.
Let us consider the following  exact sequence of homotopy groups
\begin{equation}
\pi_1(M_{k,n}) \to \pi_1 (S^1) \to \pi_0 (F) \to 0  = \pi _0 (M_{k,n}).\label{homot1}
\end{equation}

To show that $ \pi_0 (F) = 0$ it suffices  to prove that the  map $\pi_1 (M_{k,n}) \to \pi_1 (S^1)$ is surjective.  Since $p [\gamma]$ is a generator of $H^1 (M_{k,n}, \Z)$ there exists an element $a \in  H_1 (M_{k,n}, \Z)$ such that
$ \la p[\gamma] , a\ra = 1$. Since $\la [dt], \pi_* (a) \ra  =  \la  p[\gamma] , a\ra = 1$, it follows that  $\pi_*: H_1 (M_{k,n}, \Z) \to H_1 (S^1)$ is  surjective. Hence  $\pi_*: \pi_1 (M_{k,n}) \to \pi_1 (S^1)$ is 
surjective.   Hence $F$ is connected.

 
Now let $f_t$ denote  the flow on $M_{k,n}$ generated by the vector field $Z$. We note that $\Ll_Z (\gamma) = d (\gamma (Z)) = 0$, so $f_t$ respects fibration $\pi$.  Next 
we have $\Ll_Z ( \eta) = Z \rfloor n\lambda \alpha \wedge \beta + d (\eta (Z)) = 0$. Hence $f_t$   preserves also the contact form on the fiber $F$.  This proves the  first assertion. 

The second assertion
is  a consequence of Proposition \ref{solv1} and Proposition \ref{inv2}. This completes the proof of  Theorem \ref{new1}.
\end{proof}

\section* {Acknowledgement} 

H.V.L.  thanks  Alexandre Vinogradov  for explaining  the idea of their paper \cite{PV2006}  during the   conference on integrable systems  at Hradec nad Moravici in October  2010. We acknowledge  Alexandre Vinogradov kindness for   lending  H.V.L.   the  Ph.D. Thesis of  Di Pietro in  early summer 2011, which accelerated  our  work  over sections 4 and 5 considerably. H.V.L. also thanks Thiery Bouche for sending his reprint \cite{Bouche1990},  Petr Somberg for    informing her of the paper by Rumin \cite{Rumin1994}, J\"urgen Jost for  his support and his invitation to give a lecture on this subject at the Max-Planck-Institute in Leipzig in March 2011, and the ASSMS, Government College University, Lahore-Pakistan for  their hospitality and financial support during her visits in  February and September 2011, where a part of this note has been written. 
Last but not least, H.V.L.  would like to express her gratitude  to Kaoru Ono and Lorenz Schwachhofer for their helpful remarks.

\medskip

Institute of Mathematics of ASCR, Zitna 25, 11567 Praha 1, Czech Republic,  hvle@math.cas.cz,\\
Institute of Mathematics of  ASCR, Zizkova 22, 61662 Brno, Czech Republic, vanzura@ipm.cz.

\medskip



\begin{thebibliography}{99999}
\bibitem{ACFM1989}  {\sc L. C. de Andres, L. A. Cordero, M. Fernandez and J. Mencia}, Examples of four-dimensional  compact  locally conformal K\"ahler solvmanifolds, Geometriae Dedicata 29 (1989) 227-232.
\bibitem{Banyaga2007} {\sc A.  Banyaga}, Examples of non $d_\om$-exact locally conformal symplectic forms, J. Geom. 87 (2007) 1-13.
\bibitem{BK2010} {\sc G. Bande and D. Kotschick},  Contact pairs  and locally conformal  symplectic structures, arXiv.1006.0315.
\bibitem{Bouche1990}{\sc T. Bouche}, La cohomologie effective  d'une vari\'et\'e symplectique,  Bull. Sci. Math. 2.serie, 114(1990) 115-122.
\bibitem{Brylinski1988} {\sc J.C. Brylinski}, A differential complex for Poisson manifolds, JDG 28(1988), 93-114.
\bibitem{CML1998} {\sc  D. Chinea, J. Marrero, M. de Leon}, A canonical differential complex for Jacobi manifolds,
Michigan J. Math. 45(1998) 547-579.
\bibitem{PV2006} {\sc C. Di Pietro  and  A. M. Vinogradov}, A spectral sequence  associated with a  symplectic  manifold, Dokl. Akad. Nauk  413  (2007),  no. 5, 591-593, arXiv:math/0611138.
\bibitem{Pietro2006} {\sc C. Di Pietro},  Sequenza Spettrale assoiata ad una Variet\`a Simplettica, Ph.D. Thesis, Salerno, 2006.
\bibitem{Dold1972} {\sc A. Dold}, Lectures on algebraic topology,  Springer-Verlag, 1972.
\bibitem{FIL1996} {\sc M. Fernandez, R. Ibanez  and  M. de Leon}, Poisson cohomology and canonical homology of Poisson manifolds, Archivum Mathematicum, 32 (1996), 29-56.
\bibitem{FIL1998}{\sc M. Fernandez, R. Ibanez  and  M. de Leon}, Coeffective  and de Rham cohomologies of symplectic manifolds, J. of Geometry and Phys. 27 (1998), 281-296.
\bibitem{GM1988} {\sc S. I. Gelfand and Iu. I. Manin},  Methods of homological algebras,  Moscow, Nauka, 1988. 
\bibitem{GH1978} {\sc P. Griffith and  J. Harris},  Principles of algebraic geometry, John Wiley and Sohns,  1994.
\bibitem{GL1984}{\sc F. Guerida and A. Lichnerowicz}, Geometrie des algebres de Lie locales de Kirillov, J. Math. Pures et Appl. 63 (1984) 407-484.
\bibitem{HK2011}{\sc K. Hasegawa and Y. Kamishima}, Locally Conformal K\"ahler Structures on Homogeneous Spaces, arXiv:1101.3693.
\bibitem{Koszul1985} {\sc J. L. Koszul},  Crochet de Schouten-Nijenhujs et cohomologie, in ``Elie  Cartan et les Math. d'Aujour  d'Hui", Asterique hors-seire, (1985) 251-271.
\bibitem{LSV2010}{\sc H. V. Le, P. Somberg, J. Vanzura}, Poisson smooth structures on stratified symplectic spaces,
 arXiv:1011.0462.
\bibitem{Lefschetz1924} {\sc S. Lefschetz},  L'Analysis Situs et la Geometrie Algebrique, Gauthier-Villars, Paris, 1924.
\bibitem{Lepage1946} {\sc Th. N. Lepage}, Sur certaines congruences des formes alternaes, Bull. Soc. Roy. Sci., Liege 15 (1946) 21-31.
\bibitem{Lichnerowicz1977}{\sc A. Lichnerowicz}, Les Varietes de Poisson et leurs algebres de Lie associes, J. Diff.  Geometry,
12 (1977) 253-300.
\bibitem{Lychagin1979} {\sc V. V. Lychagin}, Contact geometry and non linear second-order differential equations, Russian Math. Surveys 34:1 (1979), 149-180.
\bibitem{Marle1989} {\sc C. M. Marle}, On Jacobi manifolds and Jacobi bundles, Symplectic geometry, groupoids, and integrable systems (Berkeley, CA, 1989), 227-246, Math. Sci. Res. Inst. Publ., 20, Springer, New York, 1991. 
\bibitem{McCleary2001}{\sc J. McCleary}, A user's guide to spectral sequence,  Cambridge studies in advanced mathematics, v.58,  Mathematics Springer, 2001.
\bibitem{Millionshchikov2005} {\sc D. V. Millionshchikov}, Cohomology of solvable Lie algebras and solvmanifolds, Mathematical Notes, vol. 77 (2005) 61-71.
\bibitem{Rumin1994} {\sc M. Rumin}, Formes differentielles  sur les varietes de contact, J.D.G. 39 (1994), 281-330.
\bibitem{Sawai2007} {\sc H. Sawai}, A construction of lattices on certain solvable Lie groups, Topology and its Applications 154 (2007) 3125-3134.
\bibitem{Sawai2007b} {\sc H. Sawai}, Locally conformal K\"ahler structures on compact nilmanifolds with left-invariant complex structures, Geom Dedicata 125(2007) 93-101.
\bibitem{Tischler1970} {\sc D. Tischler}, On fibering certain foliated  manifolds over $S^1$, Topology  9 (1970) 153-154.
\bibitem{TY2009} {\sc  L.S Tseng and S.T Yau}, Cohomology and Hodge Theory on Symplectic Manifolds, I, II arXiv:0909.5418, arXiv:1011.1250.
\bibitem{Vaisman1985} {\sc I. Vaisman},  Locally conformal  symplectic manifolds, Internat. J. Math.  Math. Sci.
 8 (1985) 521-536.
\bibitem{OV2010} {\sc L. Ornea and M. Verbitsky},  A report on locally conformal K\"ahler manifolds, arxiv:1002.3473.
\bibitem{Voisin2007} {\sc C. Voisin}, Hodge Theory and Complex Algebraic Geometry, I, Cambridge University Press, 2007.
\bibitem{Weil1958}{\sc A. Weil}, Introduction a l'etud\'e des vari\'et\'es  k\"ahleriennes,  Hermann, Paris, 1958.
\bibitem{Wells1986} {\sc R. Wells}, Differential analysis on complex manifolds, Springer-Verlag, 2nd edition, 1986.
\bibitem{Yan1996} {\sc D. Yan},   Hodge structures on symplectic  manifolds, Advances in Math., 120 (1996) 143-156.
\end{thebibliography}
\end{document}